\documentclass[11pt,reqno,tbtags]{article}
\usepackage{amsmath}
\usepackage{amsthm}

\numberwithin{equation}{section}
\let\sect=\section

\newtheorem{theorem}{Theorem}[section]
\newtheorem{corollary}[theorem]{Corollary}
\newtheorem{lemma}[theorem]{Lemma}
\newtheorem{proposition}[theorem]{Proposition}
\newtheorem{claim}[theorem]{Claim}
\newtheorem{example}[theorem]{\sl Example}

\theoremstyle{definition}
\newtheorem{remark}[theorem]{Remark}

\newcommand{\EE}{{\bf  E}}

\newcommand{\RR}{{\bf  R}}

\newcommand{\CC}{{\bf  C}}

\newcommand{\PP}{{\bf  P}}

\newcommand{\Var}{{\bf Var}}

\newcommand{\Lc}{{\cal L}}

\newcommand{\Leq}{{\,\stackrel{\Lc}{=}\,}}

\newcommand{\begp}{\begin{proposition}}
\newcommand{\enp}{\end{proposition}}
\newcommand{\begt}{\begin{theorem}}
\newcommand{\ent}{\end{theorem}}
\newcommand{\begl}{\begin{lemma}}
\newcommand{\enl}{\end{lemma}}
\newcommand{\begc}{\begin{corollary}}
\newcommand{\enc}{\end{corollary}}
\newcommand{\begcl}{\begin{claim}}
\newcommand{\encl}{\end{claim}}
\newcommand{\begr}{\begin{remark}}
\newcommand{\enr}{\end{remark}}
\newcommand{\begal}{\begin{algorithm}}
\newcommand{\enal}{\end{algorithm}}
\newcommand{\begd}{\begin{definition}}
\newcommand{\enf}{\end{definition}}
\newcommand{\begx}{\begin{example}}
\newcommand{\enx}{\end{example}}
\newcommand{\bega}{\begin{array}}
\newcommand{\ena}{\end{array}}

\newcommand{\sfrac}[2]{{\textstyle\frac{#1}{#2}}}

\newcommand\set[1]{\ensuremath{\{#1\}}}

\def\rompar(#1){\textup(#1\textup)}    
\newcommand\xfrac[2]{#1/#2}
\newcommand\parfrac[2]{\Bigl(\frac{#1}{#2}\Bigr)}
\newcommand\bigpar[1]{\bigl(#1\bigr)}
\newcommand\Bigpar[1]{\Bigl(#1\Bigr)}
\newcommand\biggpar[1]{\biggl(#1\biggr)}

\newcommand\expx[1]{\exp\bigl(#1\bigr)}
\newcommand\ceil[1]{\lceil#1\rceil}
\newcommand\floor[1]{\lfloor#1\rfloor}
\newcommand\ntoo{\ensuremath{n\to\infty}}
\newcommand\intoo{\int_{-\infty}^{\infty}}
\newcommand\nnoo[1]{\|#1\|_{\infty}}
\newcommand\nni[1]{\|#1\|_1}
\newcommand\eps{\varepsilon}
\newcommand\cS{\mathcal S}
\newcommand\gd{\delta}
\newcommand\gl{\lambda}
\newcommand\gs{\sigma}
\newcommand\gss{\sigma^2}

\newcommand{\refS}[1]{Section~\ref{#1}}
\newcommand{\refT}[1]{Theorem~\ref{#1}}
\newcommand{\refC}[1]{Corollary~\ref{#1}}
\newcommand{\refL}[1]{Lemma~\ref{#1}}
\newcommand{\refR}[1]{Remark~\ref{#1}}
\newcommand{\refand}[2]{\ref{#1} and~\ref{#2}}

\newcommand\ie{i.e.\spacefactor=1000}
\newcommand\eg{e.g.\spacefactor=1000}

\newcommand\cf{{cf.}\spacefactor=1000}

\newcommand\nopf{\qed}   
\newcommand\noqed{\renewcommand{\qed}{}} 
\newcommand\qedtag{\tag*{\qedsymbol}}

\newcommand\mgf{moment generating function}

\newcommand\KSm{Kolmogorov--Smirnov}
\newcommand\Holder{H\"older}
\newcommand\Roesler{R\"{o}sler}
\newcommand\Quicksort{\texttt{Quicksort}}

\newcommand\phizn{\phi_{Z_n}}
\newcommand\psizn{\psi_{Z_n}}
\newcommand\psiznl{\psizn(\gl)}
\newcommand\psiy{\psi_{Y}}
\newcommand\psiyl{\psiy(\gl)}
\newcommand\supfnf{\sup_x|f_n(x)-f(x)|}
\newcommand\dtv{d_{\mbox{\rm \scriptsize TV}}}
\newcommand\dks{d_{\mbox{\rm \scriptsize KS}}}
\begin{document}

\setcounter{page}{0}
\thispagestyle{empty}

\begin{center}
{\Large \bf Approximating the Limiting Quicksort Distribution \\}
\normalsize

\vspace{4ex}
{\sc James Allen Fill\footnotemark} \\
\vspace{.1in}
Department of Mathematical Sciences \\
\vspace{.1in}
The Johns Hopkins University \\
\vspace{.1in}
{\tt jimfill@jhu.edu} and {\tt http://www.mts.jhu.edu/\~{}fill/} \\
\vspace{.2in}
{\sc and} \\
{\sc Svante Janson}\\ 
\vspace{.1in}
Department of Mathematics \\
\vspace{.1in}
Uppsala University \\
\vspace{.1in}
{\tt svante.janson@math.uu.se} and {\tt http://www.math.uu.se/\~{}svante/} \\
\end{center}
\vspace{3ex}

\begin{center}
{\sl ABSTRACT} \\
\end{center}

The limiting distribution of the normalized number of comparisons used
by {\tt Quick\-sort} to sort an array of~$n$ numbers is known to be the unique
fixed point with zero mean of a certain distributional transformation~$S$.
We study the convergence to the limiting distribution of the sequence of
distributions obtained by iterating the transformation~$S$, beginning with a
(nearly) arbitrary starting distribution.  We demonstrate geometrically fast
convergence for various metrics and discuss some implications for
numerical calculations of the limiting {\tt Quicksort} distribution.  Finally,
we give companion lower bounds which show that the convergence is not faster
than geometric.
\bigskip
\bigskip

\begin{small}

\par\noindent
{\em AMS\/} 2000 {\em subject classifications.\/}  Primary 68W40;
secondary 68P10, 60E05, 60E10, 60F05.
\medskip
\par\noindent
{\em Key words and phrases.\/}
{\tt Quicksort}, characteristic function, density,
moment generating function,
sorting algorithm,
coupling,  
Fourier analysis,
Kolmogorv--Smirnov distance, total variation distance,
integral equation,
numerical analysis, $d_p$-metric.
\medskip
\par\noindent
\emph{Date.} January~15, 2001.
\end{small}

\footnotetext[1]{Research supported by NSF grant DMS--9803780,
and by the Acheson J.~Duncan Fund for the Advancement of Research in
Statistics.}

\newpage
\addtolength{\topmargin}{+0.5in}

\section{Introduction and summary}

The
\Quicksort{} algorithm of Hoare~\cite{Hoare} is ``one of the
fastest, the
best-known, the most generalized,
\ldots and the most
widely used algorithms for sorting an array of numbers''~\cite{ES}.
{\tt Quicksort}
is the standard sorting procedure in {\tt Unix} systems,
and in a special issue of
\emph{Computing in Science \& Engineering},
guest editors Jack Dongarra and Francis Sullivan (\cite{DS};
see also~\cite{JaJa})
chose \Quicksort{} as one of the ten algorithms
``with the greatest influence on the development and practice
of science and engineering in the 20th century."
Our goal in this introductory section is to review briefly some of
what is known about
the analysis of {\tt Quicksort} and to summarize how this paper
advances that analysis.

The {\tt Quicksort} algorithm for sorting an array of~$n$ numbers is
extremely simple
to describe.  If $n = 0$ or $n = 1$, there is nothing to do.  If $n
\geq 2$, pick a number
uniformly at random from the given array.  Compare the other numbers
to it to partition the
remaining numbers into two subarrays.  Then recursively invoke {\tt
Quicksort} on each of
the two subarrays.

Let~$X_n$ denote the (random) number of comparisons required (so that
$X_0 = 0$).  Then~$X_n$
satisfies the distributional recurrence relation
$$
X_n \Leq X_{U_n - 1} + X^*_{n - U_n} + n - 1,\qquad n \geq 1,
$$
where~$\Leq$ denotes equality in law (i.e.,\ in distribution), and
where, on the right,
$U_n$ is distributed uniformly on the set $\{1, \ldots, n\}$,
$X_j^* \Leq X_j$,
and
$$
U_n;\ X_0, \ldots, X_{n - 1};\ X^*_0, \ldots, X^*_{n - 1}
$$
are all independent.

As is well known and quite easily established, for $n \geq 0$ we have
$$
\mu_n := \EE\,X_n = 2 (n + 1) H_n - 4 n \sim 2 n \ln n,
$$
where $H_n := \sum_{k = 1}^n k^{-1}$ is the $n$th harmonic number
and~$\sim$ denotes
asymptotic equivalence.  It is also routine to compute explicitly the
standard deviation
of~$X_n$ (see Exercise 6.2.2-8 in~\cite{Knuth3}), which turns out to be
$\sim n \sqrt{7 - \frac{2}{3} \pi^2}$.

Consider the normalized variate
\begin{equation}
\label{normalized}
Y_n := (X_n - \mu_n) / n,\qquad n \geq 1.
\end{equation}
R\'{e}gnier~\cite{Reg} showed using martingale arguments that $Y_n \to Y$
in distribution, with~$Y$ satisfying the distributional identity
\begin{equation}
\label{fix}
Y\,\Leq\,U Y + (1 - U) Y^* + g(U) =: h_{Y, Y^*}(U),
\end{equation}
where
\begin{equation}\label{gu}
g(u) := 2 u \ln u + 2 (1 - u) \ln (1 - u) + 1,
\end{equation}
and where, on the right of~$\Leq$ in~\eqref{fix},
$U$, $Y$, and~$Y^*$ are independent, with $Y^* \Leq Y$  and $U \sim
\mbox{unif}(0, 1)$.  R\"{o}sler~\cite{Roesler} showed that~\eqref{fix}
characterizes the limiting law~$\Lc(Y)$, in the precise sense that
$F := \Lc(Y)$ is the \emph{unique} fixed point of the operator
\begin{equation}
\label{operator}
G = \Lc(V) \mapsto SG := \Lc(U V + (1 - U) V^* + g(U))
\end{equation}
(in what should now be obvious notation) subject to
$$
\EE\,V = 0,\qquad \Var\,V < \infty.
$$
[The fixed points of $G$ with finite mean are the translates
$\Lc(Y+c)$ with~$c$ constant, but there are other fixed points without
mean; see \cite{SJFill1A} for a complete characterization.]

R\"{o}sler~\cite{Roesler} showed that
the moment generating function
of the limiting distribution $\Lc(Y)$
is everywhere finite.
We have studied the limiting distribution further
in \cite{SJFill1},
showing
that $\Lc(Y)$
has a density~$f$ which is
infinitely differentiable, and that each derivative~$f^{(k)}(y)$
is bounded and
decays as $y \to \pm \infty$
more rapidly than any power of~$|y|^{-1}$.
(This improves an earlier result by
Tan and Hadjicostas~\cite{TanH}.)

The purpose of the present paper is to study the convergence to the
limiting distribution $\Lc(Y)$ of the sequence of distributions
obtained by iterating \Roesler's operator $S$ in \eqref{operator},
beginning with a (nearly) arbitrary starting distribution.
To fix notation, we let $Z_0$ be an arbitrary random variable, and
$F_0:=\Lc(Z_0)$ its distribution.  We define, for $n\ge1$,
$$
Z_{n}:=h_{Z_{n-1}^{\phantom{*}},Z_{n-1}^*}(U),
$$
with $Z^*_{n-1}\Leq Z_{n-1}$  and
$Z_{n-1}$, $Z^*_{n-1}$, and $U$ independent; in other words,
$$
F_n:=\Lc(Z_n)=S^n F_0,
\qquad n\ge0.
$$
Let $\|X\|_2:=(\EE\,X^2)^{1/2}$ denote the $L^2$-norm, and let
$d_2$ denote the metric on the space of probability distributions
with finite variance defined by
\begin{equation}\label{d2}
d_2(F,G):=\min\|X-Y\|_2,
\end{equation}
taking the minimum over all pairs of random variables $X$ and $Y$
(defined on the same probability space) with $\Lc(X)=F$ and
$\Lc(Y)=G$.
Note that, using the coupling with $X$ and $Y$ independent,
for any~$F$ and~$G$ each with zero mean and finite variance,
\begin{equation}\label{d2est}
d_2(F,G)
\le (\EE\,X^2 + \EE\,Y^2)^{1/2}
\le\|X\|_2+\|Y\|_2,
\end{equation}
when~$\Lc(X)=F$ and $\Lc(Y)=G$.
\Roesler{}~\cite{Roesler} showed that if $Z_0$ has mean 0 and
finite variance, then
$F_n\to F$ in the $d_2$-distance with a geometric rate:
\begin{equation}\label{d2conv}
d_2(F_n,F) \le (2/3)^{n/2} d_2(F_0,F)
\le (2/3)^{n/2}(\Var\,Z_0 + \gss)^{1/2},
\end{equation}
where
\begin{equation}\label{gss}
\sigma^2:=\Var\,Y = 7 - \sfrac{2}{3} \pi^2 \doteq 0.42.
\end{equation}
Our main interest is to show similar estimates for other measures of
the distance between $F_n$ and $F$.

We will show in \refS{S:convdens},
using estimates of the characteristic functions given in
\cite{SJFill1} and \refS{S:chf},
that the distribution $F_n$ has a bounded, continuous density function
$f_n$, at least as soon as $n\ge3$, and that, if $Z_0$ has mean
0 and finite variance, then $f_n$ converges uniformly
to $f$, with a geometric rate of convergence, as \ntoo.
We further show geometrically fast convergence in the total variation
and \KSm{} distances, too.

In \refS{S:mgf1} we give bounds for the moment generating functions of
$Y$ and of $Z_n$.
In \refS{S:mgf2} we show that if $Z_0$ has mean 0 and a finite
\mgf{}~$\psi_0$,
then the \mgf{}~$\psi_n$ of $F_n$ is finite and
converges uniformly on compact intervals to the \mgf{} of~$Y$, again with
a geometric rate of convergence.
We study in particular the cases $Z_0 = 0$
and~$Z_0$ normally distributed with zero mean and sufficiently large variance;
it turns out that in these cases~$\psi_n(\lambda)$ converges monotonically.

In \refS{S:numerics}, we discuss some implications for
numerical calculations of the limiting \Quicksort{} distribution $F$,
showing how explicit and arbitrarily small error bounds can be obtained.

Finally, in Sections~\ref{S:lower}--\ref{S:otherlower} we give some
companion lower bounds, showing that the convergence is not faster
than geometrical for several different metrics.
We also show geometrically fast convergence in the $d_p$ metric for
any finite $p$.

\begin{remark}
The mode and rate of convergence of the distribution of the actual normalized
\Quicksort{} variables~$Y_n$ of~(\ref{normalized}) to
the limit $F$ is a quite different matter, which will be studied in
another paper \cite{SJFill3}.
\end{remark}

\section{Bounds on the characteristic functions}\label{S:chf}
In \cite{SJFill1} we gave bounds on the characteristic function of
$Y$.
The same method yields, more generally, bounds on the characteristic
function of $Z_n$ for arbitrary $Z_0$.
We write $\phi_X(t):=\EE\,e^{itX}$ for any random variable~$X$.

\begin{theorem}\label{T:phi}
For every real $p \geq 0$ there is a constant $0 < c_p < \infty$ such
that for any $Z_0$ and any $n>p+1$,
the characteristic function
$\phi_{Z_n}(t)$ satisfies
\begin{equation}\label{phip}
|\phi_{Z_n}(t)|\le c_p|t|^{-p}\mbox{\rm \ \ for all $t \in \RR$.}
\end{equation}
The best possible constants $c_p$ satisfy
$c_0=1$,
$c_{1/2}\le2$,
$c_{3/4}\le\sqrt{8\pi}$,
$c_1\le4\pi$,
$c_{3/2}<187$, $c_{5/2}<103215$,
$c_{7/2}<197102280$,
and the relation
\begin{equation} \label{cp+1}
c_{p+1}\le 2^{p+1} c_p^{1+(1/p)} p/(p-1),
 \qquad p>1;
\end{equation}
moreover, at least if we restrict \eqref{phip} to $n\ge p+2$,
\begin{equation} \label{cp}
c_p\le 2^{p^2+6p},\qquad p>0.
\end{equation}
\end{theorem}
[The bounds on the constants $c_p$ obtained here are the same as for
the special case $Z_0 \Leq Y$ (whence $Z_n \Leq Y$ for every~$n$)
in~\cite{SJFill1}. However, there is no reason to believe that our method
yields the best possible bounds, and the best constants
for the special case in~\cite{SJFill1}
may be smaller than the best constants in \refT{T:phi} here.]

\begin{proof}
The proof is almost identical to the proof of the special case in
\cite{SJFill1}, so we will omit some details.
For any random variable~$Z$,  we abuse notation slightly
and denote by~$S Z$ the random variable
$h_{Z,Z^*}(U)=U Z + (1 - U) Z^* + g(U)$
where
$U$, $Z$, and~$Z^*$ are independent, with $Z^* \Leq Z$  and $U \sim
\mbox{unif}(0, 1)$; thus $S Z$ is a random variable with the
distribution $S \Lc(Z)$.
By conditioning on $U$, we obtain the fundamental relation
\begin{equation}
\label{phiinteq}
\phi_{S Z}(t)
= \int^1_0\!\phi_Z(u t)\,\phi_Z((1 - u) t)\,e^{i t g(u)}\,du,
\qquad t \in \RR,
\end{equation}
and thus the estimate
\begin{equation}
\label{phirecbd}
|\phi_{S Z}(t)|
\leq \int^1_0\!|\phi_Z(u t)|\,\,|\phi_Z((1 - u) t)|\,\,du.
\end{equation}

To complete the proof, we give a series of lemmas.

\begin{lemma}\label{L:hyz}
For any real numbers~$y$ and~$z$, the random variable~$h_{y,z}(U)$
defined by \eqref{fix} satisfies
\begin{equation*}
|\EE\,e^{ith_{y,z}(U)}|\le 2|t|^{-1/2}.
\end{equation*}
\end{lemma}

\begin{proof}
This follows by a method of van der Corput  \cite{vdC,Monty,SJFill1},
using little more than the fact that $h_{y, z}$ is convex with
$h''_{y, z} \geq 8$ on $(0,1)$.
\end{proof}

\begin{lemma}\label{L:1/2}
For any random variable $Z$ and real $t$, we have
$|\phi_{S Z}(t)|\le 2|t|^{-1/2}$.
\end{lemma}

\begin{proof}
\refL{L:hyz} yields
$$
|\phi_{S Z}(t)|=
\Bigl|\EE\,e^{ith_{Z,Z^*}(U)}\Bigr|
\le \EE\Bigl|\EE\Bigl(e^{ith_{Z,Z^*}(U)}\Bigm| Z,Z^*\Bigr)\Bigr|
\le 2 |t|^{-1/2}.
\eqno\qedsymbol
$$
\noqed
\end{proof}

Returning to our sequence $(Z_n)$,
the preceding lemma applies to all elements except $Z_0$,
\ie{},
\begin{equation}\label{zn1/2}
|\phizn(t)|\le 2|t|^{-1/2},
\qquad n\ge1,
\end{equation}
which yields the case $p = 1/2$ of~\refT{T:phi}.
We improve the exponent by induction, using~\eqref{phirecbd}.

\begin{lemma}\label{L:p01}
Let $0<p<1$.
If $|\phi_Z(t)|\le c_p |t|^{-p}$, $t \in \RR$, then
$$
|\phi_{S Z}(t)|
\le \frac{\bigl[\Gamma(1-p)\bigr]^2}{\Gamma(2-2p)}  c_p^2 |t|^{-2p}.
$$
\end{lemma}

\begin{proof}
By~\eqref{phirecbd} and the hypothesis,
\begin{equation*}
|\phi_{S Z}(t)|
\le
\int_0^1\!c_p^2|ut|^{-p}|(1-u)t|^{-p}\,du
=c_p^2 |t|^{-2p}\int_0^1\!u^{-p}(1-u)^{-p}\,du,
\end{equation*}
and the result follows by evaluating the beta integral.
\end{proof}

In particular, using \eqref{zn1/2},
\refL{L:p01} yields
\begin{equation}
\label{secbd}
|\phizn(t)| \leq 4 \pi|t|^{-1},
\qquad n\ge2.
\end{equation}
This proves~\eqref{phip} for $p = 1$, with $c_1 \leq 4 \pi$.
Since $|\phizn(t)|\le1$, for any $p\le1$ we trivially have
$|\phizn(t)|\le|\phizn(t)|^p$, which by~\eqref{secbd}
establishes~\eqref{phip} for all
$p \le 1$ with $c_p \leq (4 \pi)^p$;
applying \refL{L:p01} again, we obtain \eqref{phip} for all $p < 2$.
Somewhat better numerical bounds are obtained for $1/2 < p < 1$
by taking a geometric average between the cases $p = 1/2$ and $p = 1$;
this yields
$c_p \leq 2^{2p} \pi^{2p-1}$, $1/2 \le p \le 1$.
In particular, we have $c_{3/4} \leq \sqrt{8\pi}$, and thus, by
\refL{L:p01},
$c_{3/2} \leq 8\pi^{1/2} \bigl[\Gamma(1/4)\bigr]^2 < 186.4 < 187$.

\begin{lemma}\label{L:p>1}
Let $p > 1$.
If $|\phi_Z(t)|\le c_p |t|^{-p}$, $t \in \RR$, then
$$
|\phi_{S Z}(t)|
\le 2^{p+1} c_p^{1 + (1/p)}\frac p {p - 1}|t|^{-(p+1)}.
$$
\end{lemma}

\begin{proof}
This is similar to the proof of \refL{L:p01}, substituting the
hypothesis
(and the trivial $|\phi_Z|\le1$) into \eqref{phirecbd}, but the estimate
of the integral
is slightly more complicated;
for details see \cite{SJFill1}.
\end{proof}

\refL{L:p>1} completes, by induction,
the proof of \eqref{phip}
and the estimate~\eqref{cp+1}.

The bound for~$c_{3/2}$ obtained above and \eqref{cp+1} now yield
(using {\tt Maple}) first
$c_{5/2} < 103215$ and then $c_{7/2} < 197102280$.
These bounds and \eqref{cp+1} further yield
\begin{equation} \label{cp3/2+}
c_p\le 2^{p^2+5p},
\qquad p = k + \sfrac{3}{2},
\end{equation}
for integers $k\ge0$;
again see \cite{SJFill1} for details.
To obtain \eqref{cp}
if $p>1/2$, let $p_1:=\ceil{p-\frac12}+\frac12$. Then,
by \eqref{phip} and \eqref{cp3/2+}, provided $n\ge p+2>p_1+1$,
$$
|\phizn(t)|^{1/p}
\le |\phizn(t)|^{1/p_1}
\le 2^{p_1+5}|t|^{-1}
\le 2^{p+6}|t|^{-1}.
$$
The case $p\le1/2$ follows similarly from \eqref{zn1/2},
which completes the proof of \refT{T:phi}.
\end{proof}

\begin{remark}
A variety of other bounds are possible.  For example, if we begin with the
inequality~\eqref{secbd} and use~\eqref{phirecbd},
we can easily derive the following result:
\begin{equation}
\label{thirdbd}
|\phizn(t)|
\leq \frac{32 \pi^2}{t^2} \Bigl(\ln \Bigl(\frac{t}{4\pi}\Bigr)+2\Bigr)
\leq \frac{32 \pi^2 \ln t}{t^2}\mbox{\ \ for all $t \geq 1.72$ and $n\ge3$.}
\end{equation}
\end{remark}

\section{Convergence of densities} \label{S:convdens}

It is easily checked that the random variable $h_{y,z}(U)$ is
absolutely continuous for
every fixed $y$ and $z$, and thus, by mixing, $S Z$ is absolutely
continuous for every $Z$.
In other words, for any $Z_0$, the random variables $Z_n$ have
densities for all $n\ge1$; \cf{} \cite{TanH}.
These densities may be unbounded and discontinuous, at least for
$n=1$, as is seen in the case $Z_0\equiv0$. However, we now can show
that for $n\ge3$, at least, no such irregularities occur.

\begin{theorem}\label{T:density}
If $n\ge3$, then $Z_n$ has a bounded continuous density function
$f_n$, for any $Z_0$.
More generally, if $k\ge0$, then $f_n$ is $k$ times continuously
differentiable for all $n\ge k+3$, and there exists a constant $C_k$
independent of $Z_0$ and $n$ \rompar(with $n\ge k+3$) such that
$|f_n^{(k)}(x)|\le C_k$, $x\in\RR$.
Explicitly, $|f_n(x)|\le 16$ when $n\ge5$, and
$|f_n'(x)|\le 2466$ when $n\ge6$.
\end{theorem}

\begin{proof}
\refT{T:phi} shows, in particular, that as soon as $n\ge3$,
$$
|\phizn(t)|\le \min(1,187|t|^{-3/2}),
$$
and thus $\phizn$ is
integrable.
This implies, as is well-known (see \eg, \cite[Theorem~XV.3.3]{Feller2})
that~$Z_n$ has a bounded continuous density~$f_n$ given by the
Fourier inversion formula
\begin{equation}
\label{Fourierinv}
f_n(x) = \frac{1}{2 \pi}\int^{\infty}_{-\infty}\!e^{- i tx}\,\phizn(t)\,dt,
\qquad x \in \RR.
\end{equation}
Moreover, using \refT{T:phi} with $p=k + \frac{3}{2}$, we see that
$t^k\phizn(t)$ is also integrable when $n\ge k+3$, which by a standard argument
shows
that~$f_n$ is $k$ times differentiable, with
\begin{equation}\label{Fourierinvk}
f_n^{(k)}(x)
= \frac{1}{2 \pi}\int^{\infty}_{-\infty}\!
(- i t)^k\,e^{- i t x}\,\phizn(t)\,dt,
\qquad x \in \RR;
\end{equation}
and thus
\begin{equation}
\label{fkbd}
\sup_x |f_n^{(k)}(x)|
\le \frac1{2 \pi} \int_{-\infty}^\infty\!|t|^k\,|\phizn(t)|\,dt,
\end{equation}
where the latter integral can be estimated using \refT{T:phi} with $p=k +
\frac{3}{2}$.

The argument above
yields the bound
\begin{equation}\label{fn3/2}
|f_n(x)|
\le \frac1{2 \pi} \int_{-\infty}^\infty\!\min(1,187|t|^{-3/2})\,dt
=\frac3{\pi} 187^{2/3}<31.3,
\qquad n\ge3.
\end{equation}
To obtain better numerical bounds we combine
\refT{T:phi} for $p=0$, $1/2$, $3/2$, $1$, $5/2$, $7/2$ and~\eqref{thirdbd}
(for $t$ in different intervals; see  \cite{SJFill1} for details);
this yields, provided $n \ge 5$,
$f_n(x)\le \frac1{2\pi}\int\!|\phi_n| <15.3$; similarly, invoking also
\eqref{phip} with $p=9/2$,
$f'_n(x)\le \frac1{2\pi}\int\!|t||\phi_n(t)|\,dt <2465.9$ for $n\ge6$.
\end{proof}

\begin{theorem}\label{T:densconv}
Suppose that $\EE\,Z_0=0$ and $\Var\,Z_0<\infty$.
Then the density functions~$f_n$ of  \refT{T:density}
converge uniformly to
the (smooth) density function~$f$ of~$Y$ at a geometric rate:
$$
\supfnf = O(r^n)\ \ \mbox{\rm for every fixed $r > (2/3)^{1/2}$.}
$$
Explicitly, for any $p>1$ and $n>p+1$,
\begin{equation}\label{fn1}
\sup_x|f_n(x)-f(x)|
\le \frac{A}{2\pi} \parfrac{2c_p}{A}^{2/(p+1)}\frac{p+1}{p-1}
  \parfrac23^{(\frac12-\frac1{p+1})n},
\end{equation}
where $A:=(\Var\,Z_0+\sigma^2)^{1/2}$ and $c_p$ is as in \refT{T:phi}.
In particular,
\begin{equation}\label{fn2}
\sup_x|f_n(x)-f(x)|
\le 2297 A \parfrac23^{5n/18}
<  2297 A (0.8935)^n,
\qquad n\ge5.
\end{equation}
Moreover,
\begin{equation}\label{fn3}
\sup_x|f_n(x)-f(x)|
\le \frac{128A}{\pi} \parfrac23^{(n / 2) - 3.7 \sqrt n},
\qquad n\ge3.
\end{equation}
\end{theorem}

\begin{proof}
By the Fourier inversion formula \eqref{Fourierinv},
\begin{equation}\label{diff1}
|f_n(x) - f(x)| \leq
\frac{1}{2 \pi} \int^{\infty}_{-\infty}\!|\phizn(t) - \phi_Y(t)|\,dt.
\end{equation}
In order to estimate the right hand side, note that
for any random variables $X$ and $Y$,
$$
|\phi_X(t)-\phi_Y(t)|\le  \EE|e^{itX}-e^{itY}|\le \EE|tX-tY|
\le |t|\,\|X-Y\|_2;
$$
since the characteristic functions here depend on the marginal
distributions only,
this and the definition~\eqref{d2} yield
$$
|\phi_X(t)-\phi_Y(t)|
\le |t| d_2(\Lc(X),\Lc(Y)).
$$
In particular, with $d_n := d_2(F_n, F)$,
\begin{equation}
\label{julie}
|\phizn(t)-\phi_Y(t)|
\le |t| d_n.
\end{equation}
Further, for any $p>1$ and $n>p+1$, \refT{T:phi} yields the estimate
$$
|\phizn(t)-\phi_Y(t)|
\le |\phizn(t)|+|\phi_Y(t)|
\le 2c_p|t|^{-p}.
$$
Consequently, for any $ T>0$,
\begin{align*}
\int^{\infty}_{-\infty}\!|\phizn(t) - \phi_Y(t)|\,dt
 &\leq \int^T_{-T}\!d_n|t|\,dt
  +  \int_{|t| > T}\!2c_p |t|^{- p}\,dt \\
 &  =  d_n T^2 + 4 \frac{c_p}{p - 1} T^{1- p}.
\end{align*}
For given~$n$ and~$p$, the optimal choice here is
$
T:= \left( \xfrac{2 c_p}{d_n} \right)^{1 / (p + 1)}
$,
giving the bound
\begin{equation}\label{axel}
\int^{\infty}_{-\infty}\!|\phizn(t) - \phi_Y(t)|\,dt
\le
\frac{p+1}{p-1}(2c_p)^{2/(p+1)}  d_n^{1 - (2/(p+1))}.
\end{equation}
With \eqref{diff1} and the estimate \eqref{d2conv}, this yields \eqref{fn1}.
Choosing $p=7/2$ and evaluating the constants numerically, using
$A\ge \sigma>0.648$, we obtain \eqref{fn2}.

To obtain the final estimate, we use \eqref{cp}
and observe that, for $p\ge2$,
$$
\parfrac{2c_p}{A}^{2/(p+1)} \le
\parfrac{2^{p^2+6p+1}}{\sigma}^{2/(p+1)}
=2^{2(p+5)} \parfrac{2^{-4}}{\sigma}^{2/(p+1)} \le
2^{2(p+5)} \Bigpar{1-\frac2{p+1}},
$$
which by \eqref{fn1} yields that for $n\ge p+2 \ge4$,
\begin{equation*}
\sup_x|f_n(x)-f(x)|
\le \frac{A}{2\pi} 2^{2p+10}
  \parfrac23^{(\frac12-\frac1{p+1})n}.
\end{equation*}
Choosing the optimal $p:=\left[ n \ln(3/2) / (2 \ln2) \right]^{1/2} - 1$,
we find~\eqref{fn3} [with the constant
$(8 (\ln2) / \ln(3/2))^{1/2}<3.69812$ multiplying $\sqrt n$],
at least when $n\ge31$. For  $3\le n\le30$, \eqref{fn3} follows trivially
from
\eqref{fn3/2}, since the right hand side of
\eqref{fn3} then is larger than $193$.
\end{proof}

To test out
\refT{T:densconv} numerically, choose $Z_0 \equiv 0$, so that
$A = \sigma \doteq 0.648$.  For $n = 100$, \eqref{fn2} yields the bound
$0.0192$;
for $n \geq 177$, \eqref{fn3} is better, and yields for example
$3.21 \times 10^{-6}$ for $n = 177$, $2.07 \times 10^{-6}$ for $n = 180$,
and $1.07 \times 10^{-7}$ for $n = 200$.

\begin{remark}
Similarly, using \eqref{Fourierinvk}, we obtain geometric uniform
convergence of the first derivatives, and of any higher derivatives,
of the density functions.
\end{remark}

\begin{remark}
Suppose that~$Z_0$ has finite moments of all orders.
Then, by \refL{L:decrease} below, $\EE |Z_n|^p$ is finite and stays
bounded in~$n$, for each
real $0 \leq p < \infty$.  It
follows that the characteristic functions~$\phi_{Z_n}$ are infinitely
differentiable
with derivatives bounded uniformly in~$n$.
If we apply both \refT{T:phi} and~\eqref{julie} to $|\phi_{Z_n}(t) -
\phi_Y(t)|$ and take
the geometric mean of the resulting bounds, we find, for
$n > 2 p + 2$,
$$
\left| \phi_{Z_n}(t)-\phi_Y(t) \right|
\leq \left[ 2c_{2p+1} |t|^{-2p} d_2(F_n,F)\right]^{1/2}.
$$
It follows easily by induction on $k$, using \cite[Lemma 2.10]{SJFill1},
that in fact, for every real $p \geq 0$ and
integer $k \geq 0$, there is a constant $c_{p, k}$ [depending on
$\Lc(Z_0)$] such that for
all
$n > 2^{k+1}p + 2$ we have, with $\rho_k := (2/3)^{2^{-k-2}}<1$,
$$
\sup_{t \in \RR} |t|^p\left|\phi^{(k)}_{Z_n}(t)-\phi^{(k)}_Y(t)\right|
\leq c_{p, k}\rho_k^n.
$$
Omitting details,
since the Fourier transform is continuous on the Schwartz space
\cite{Schwartz}
$$
\cS := \set{f:\sup_t |t|^p |f^{(k)}(t)|<\infty\mbox{\ for all\ }p,k\ge0},
$$
it follows that for each $k$ and $p$,
$|x|^p f^{(k)}_n(x)$ converges uniformly to $|x|^pf^{(k)}(x)$ with
geometric rate.
\end{remark}

\refT{T:densconv} treats uniform approximation of $f$ by $f_n$, using
the norm $\nnoo{f_n-f}:=\supfnf$. We now turn to studying the
error in the $L^1$-norm
$\nni{f_n-f}:=\intoo\!|f_n-f|$.

Note first that, because $\intoo\!f_n=\intoo\!f=1$,
\begin{equation*}
\tfrac12\intoo\!|f_n(x)-f(x)|\,dx
=\intoo\!\bigpar{f(x)-f_n(x)}^+\,dx,
\end{equation*}
and that this coincides with the total variation distance
\begin{equation*}
\dtv(F_n,F):=\sup_{A \subseteq \RR}|\PP(Z_n\in A)-\PP(Y\in A)|;
\end{equation*}
moreover, it dominates the \KSm{} distance
\begin{equation*}
\dks(F_n,F):=\sup_{ x\in\RR}|\PP(Z_n\le x)-\PP(Y\le x)|
\le\dtv(F_n,F).
\end{equation*}

\begin{theorem}\label{T:dtvconv}
Suppose that $\EE\,Z_0 = 0$ and $\Var\,Z_0<\infty$.
Then the total variation and \KSm{} distances between $F_n$ and $F$
converge geometrically to $0$:
$\dks(F_n,F)\le\dtv(F_n,F)=O(r^n)$ for every fixed $r>(2/3)^{1/2}$.
Explicitly, for any $n\ge1$,
\begin{equation}\label{fn1a}
\dks(F_n,F)\le\dtv(F_n,F)
\le 135An \parfrac23^{(n / 2) - 3.7 \sqrt n}.
\end{equation}
\end{theorem}

\begin{proof}
For any $a\in(0,1)$,
\begin{equation}\label{g1}
\dtv(F_n,F)=\intoo\!\bigpar{f(x)-f_n(x)}^+\,dx
\le \nnoo{f_n-f}^{1-a}\intoo\!f(x)^a\,dx,
\end{equation}
where $\nnoo{f_n-f}$ is estimated in \refT{T:densconv}.
The final integral can be estimated by \Holder's inequality:
for any $b>0$
\begin{equation}\label{fab}
\begin{split}
\intoo\!f(x)^a\,dx
&=\intoo\!f(x)^ae^{ab|x|}\cdot e^{-ab|x|}\,dx\\
&\le \biggpar{\intoo\!f(x)e^{b|x|}\,dx}^a
    \biggpar{\intoo\!e^{-ab|x|/(1-a)}\,dx}^{1-a}\\
&\le \left[ \psi(b)+\psi(-b) \right]^a \parfrac{2(1-a)}{ab}^{1-a}\\
&= \frac2{ab} \Bigpar{ab\frac{\psi(b)+\psi(-b)}2}^a (1-a)^{1-a},
\end{split}
\end{equation}
where $\psi(\lambda):=\EE\,e^{\lambda Y}$ is the \mgf{} of~$Y$.
\Roesler{} \cite{Roesler} proved that $\psi(\lambda)$ is finite for all
$\lambda$;
thus $\int\!f^a<\infty$ for every $a \in (0, 1)$, and the first claim
follows by
\eqref{g1} and \refT{T:densconv}.

For \eqref{fn1a} we choose $b=1/3$,
for which it will be shown in \refT{T:KL} below that
$\psi(\pm b)\le \exp(1/9)<1.2$,
and thus \eqref{fab} implies $\intoo\!f^a < 2 / (a b) = 6/ a$.
Denoting the right hand side of \eqref{fn3} by $B$, we thus obtain
from \eqref{g1} and \eqref{fn3}, observing that $B\ge (3/2)^{-n/2}$,
$$
\dtv(F_n,F)
\le \frac{6}a B^{1-a}
\le \frac{6}a (3/2)^{an/2} B.
$$
We optimize by taking $a:=2 / (n \ln(3/2))$ and obtain the following bound
(for $n\ge5$, so that $a<1$;
smaller $n$ are trivial since $\dtv\le1$):
$$
\dtv(F_n,F)
\le 3 e n B \ln(3/2)
= \frac{384 e\ln(3/2)}{\pi} A n  \parfrac23^{(n / 2) - 3.7\sqrt n}.
\eqno\qedsymbol
$$
\noqed
\end{proof}

\begin{remark}
If we are content with a weaker explicit bound, we can avoid invoking
estimates of $\psi$ by using moments of $Y$ instead. For example,
\begin{equation*}
\begin{split}
\intoo\!f(x)^{1/2}\,dx
&\le \biggpar{\intoo\!f(x)(\gss+x^2)\,dx}^{1/2}
    \biggpar{\intoo\!\frac{dx}{\gss+x^2}}^{1/2}
=(2\pi\gs)^{1/2} <2.1
\end{split}
\end{equation*}
and thus
$$
\dks(F_n,F)\le\dtv(F_n,F)\le 2.1 \nnoo{f_n-f}^{1/2}.
$$
\end{remark}

\sect{Bounds on moment generating functions}
\label{S:mgf1}
Letting $\psi_Z(\lambda):=\EE\,e^{\lambda Z}$ denote the \mgf{} of a
random variable
$Z$,
we find in analogy with \eqref{phiinteq} the relation
\begin{equation}
\label{psiinteq}
\psi_{S Z }(\lambda)
= \int^1_0\!\psi_Z(u \lambda)\,\psi_Z((1 - u) \lambda)\,e^{ \lambda g(u)}\,du,
\qquad \lambda \in \RR.
\end{equation}
In particular, it follows that if $\psi_Z(\lambda)$ is finite for all
$\lambda$,
then so is $\psi_{S Z }(\lambda)$.

\Roesler{} \cite{Roesler} proved that the \mgf{}~$\psi_Y$ is everywhere
finite and that for every $L \geq 0$ there is a
constant $K_L$ such that
\begin{equation}\label{psi1}
\psi_Y(\gl)\le e^{K_L\gl^2},
\qquad |\gl|\le L.
\end{equation}
Moreover, it is implicit in the proof that
\begin{equation}\label{psi2}
\text{if $\psi_Z(\gl)\le e^{K_L\gl^2}$ for  $|\gl|\le L$,
then $\psi_{S Z}(\gl)\le e^{K_L\gl^2}$ for  $|\gl|\le L$.}
\end{equation}
Note
that \eqref{psi2} implies by induction that if we choose  $Z_0
\equiv 0$, then
$\psi_{Z_n}(\gl)\le e^{K_L\gl^2}$, $|\gl|\le L$, for every $n$, and
thus \eqref{psi1}
follows by Fatou's lemma.
More generally, if \eqref{psi2} holds and
$\psi_{Z_0}(\gl)\le e^{K_L\gl^2}$, $|\gl|\le L$, then by induction
$\psi_{Z_n}(\gl)\le e^{K_L\gl^2}$, $|\gl|\le L$, for every $n$.

\Roesler{} did not give explicit values of the constants $K_L$,
but such values can be obtained from his proof as follows.
[Actually,  \Roesler{} \cite{Roesler} treated the somewhat more
complicated case of the
variables
$Y_n$ of~(\ref{normalized}); see \cite{SJFill3} for explicit constants
in that case.
In our case there are some simplifications leading to better constants.
Moreover, we introduce some deviations from \Roesler's proof
designed
to improve our bounds.]
\begin{theorem}
\label{T:KL}
Let $L_0\doteq 5.018$ be the largest root of $e^L=6 L^2$. Then
\eqref{psi1} and \eqref{psi2} hold with
\begin{equation*}
K_L =
\begin{cases}
1,            & L \le 0.42,\\
12,           & 0.42 < L \le L_0,\\
2 L^{-2} e^L, & L_0 < L,
\end{cases}
\end{equation*}
or any larger number. In particular, we can always take
$K_L = \max(2 L^{-2} e^L, 12)$.
\end{theorem}

For $\gl\le0$, we can obtain much better estimates.
[For~\eqref{psi2}, we restrict to $\gl\le0$ in both the assumption and
the conclusion.]
\begin{theorem}
\label{T:KL_neg}
We have
\eqref{psi1} and \eqref{psi2} for $\lambda \leq 0$ with
\begin{equation*}
K_L =
\begin{cases}
0.5,            & L \le 0.62,\\
1.25,            & 0.62 < L,
\end{cases}
\end{equation*}
or any larger number. In particular, we can always take
$K_L = 1.25$ for $\gl\le0$.
\end{theorem}

\begin{proof}[Proof of Theorems \refand{T:KL}{T:KL_neg}]
If  $\psi_Z(\gl)\le e^{K\gl^2}$ for  $|\gl|\le L$,
then by \eqref{psiinteq}, for  $|\gl|\le L$,
$$
\psi_{S Z}(\gl)\le \int_0^1 e^{K\gl^2[u^2+(1-u)^2]+\gl g(u)}\,du
=e^{K\gl^2}
\EE e^{\gl g(U)-2K\gl^2U(1-U)}.
$$
Hence, \eqref{psi2} holds with $K_L=K$
if (and only if)
\begin{equation}\label{fkl}
f_K(\gl):=\EE e^{\gl g(U)-2K\gl^2U(1-U)}
\le1,\qquad\text{when }|\gl|\le L.
\end{equation}
Similarly,
\eqref{psi2} holds with $K_L=K$ for $\gl\ge0$
(respectively, for $\gl\le0$)
if \eqref{fkl} holds for
$0\le\gl\le L$ (resp.,\ for $-L \le\gl\le 0$).
Clearly, $f_K(\gl)$ decreases as $K$ increases, and thus if some $K$
satisfies \eqref{fkl}, then so does any larger $K$.

Following \Roesler, we argue differently for small and large $L$ in
order to find a $K$ satisfying \eqref{fkl}.
For small $L$ we use a Taylor expansion. By straightforward
differentiations,
\begin{align*}
f_K(0)&=1,\\
f_K'(0)&=\EE g(U)=0,\\
f_K''(0)&=\EE\bigpar{g(U)^2-4KU(1-U)}=\tfrac13\gss-\tfrac23 K,\\
f_K'''(\gl)&=\EE\Bigl[\Bigpar{\bigpar{g(U)-4K\gl U(1-U)}^3-12KU(1-U)
  \bigpar{g(U)-4K\gl U(1-U)}}\\
&\qquad\qquad\qquad\qquad\times\expx{\gl g(U)-2K\gl^2 U(1-U)}\Bigr].
\end{align*}
We write the last formula as
$
f_K'''(\lambda) = \EE[X(U, \lambda)]
$
and note that
$0\le U (1 - U) \leq 1 / 4$ and $- \eta \leq g(U) \leq 1$,
where
$$
\eta := -g(\tfrac12)=2 \ln 2 - 1 \doteq 0.386.
$$

Consider first $\gl\ge0$.
By Taylor's formula, for $0\le\gl\le L$,
\begin{equation*}
\begin{split}
f_K(\gl)&\le
1+\tfrac12\gl^2 f_K''(0)
+\tfrac16\gl^3\sup_{0\le\gl\le L} f_K'''(\gl)\\
&\le 1+\tfrac16\gl^2\bigpar{\gss-2K+L\sup_{0\le\gl\le L} f_K'''(\gl)}
\end{split}
\end{equation*}
so \eqref{fkl} is satisfied for $\gl\ge0$ provided
\begin{equation}
\label{tomas}
L\sup_{0\le\gl\le L} f_K'''(\gl)
\le 2K-\gss.
\end{equation}

If $g(U) \geq 0$, we find
$$
X(U, \lambda) \leq (1 + 3 K^2 L) e^L, \qquad 0 \leq \lambda \leq L;
$$
while if $g(U) \leq 0$, we find
$$
X(U, \lambda) \leq 3 K ( \eta+ KL),  \qquad 0 \leq \lambda \leq L.
$$
For $K\ge1$, in either case, because $3\eta>1$,
$$
X(U, \lambda) \leq (3K\eta+3K^2L) e^L, \qquad 0 \leq \lambda \leq L,
$$
and thus
$$
L\sup_{0\le\gl\le L} f_K'''(\gl)\leq L(3K\eta+3K^2L) e^L.
$$
It is readily checked that this is less than $2K-\gss$
so that \eqref{tomas} holds, for $K=1$ and $L = 0.42$.

For larger $L$, we begin by another crude estimate. Let $W\Leq U/2$ be
uniformly
distributed on  $(0,1/2)$.
Then, by $|g(U)|\le1$ and symmetry,
\begin{equation} \label{gkl}
\begin{split}
f_K(\gl)&
\le e^{|\gl|}\EE \expx{-2K\gl^2U(1-U)}
=  e^{|\gl|}\EE \expx{-2K\gl^2W(1-W)}\\
& \le e^{|\gl|}\EE \expx{-K\gl^2W}
=e^{|\gl|}\int_0^1\expx{-K\gl^2u/2}\,du\\
&=e^{|\gl|}\frac{1-\expx{-K\gl^2/2}}{K\gl^2/2}
=: g_K(\gl).
\end{split}
\end{equation}
Note that that $g_K$, too, decreases if $K$
is increased.
Taking the logarithmic derivative, we find for $\gl>0$,
\begin{equation} \label{g'}
\begin{split}
\bigpar{\ln g_K(\gl)}'
&=1-\frac2\gl+K\gl e^{-K\gl^2/2}
\bigpar{1-\expx{-K\gl^2/2}}^{-1}\\
&=1-\frac2\gl+\frac{K\gl}
{ e^{K\gl^2/2}-1}.
\end{split}
\end{equation}
For $\gl\ge2$, this is evidently positive, and thus
$g_K$ then is increasing. Hence,
if $K\ge \tilde{K}:=2L^{-2}e^L$,
then
$$
g_K(\gl) \le g_K(L) \le g_{\tilde{K}}(L)=1-\exp(-e^{L})<1,
\qquad  2\le \gl\le L.
$$

For smaller $\gl$, we take $K = 12$, and check numerically that
$g_{12}(0.42) < 1$.  Moreover,
$$
e^{K \gl^2 / 2} - 1 = e^{6 \gl^2} - 1 \ge 6 \gl^2 + 18 \gl^4
$$
and further, if $1/3 \le \gl \le 1$,
$$
(1-\sfrac{\gl}{2})^{-1} \le 1 + \gl \le 1 + 3 \gl^2
$$
and thus
$$
\frac{K \gl}{e^{K \gl^2 / 2} - 1} \le \frac{12 \gl}{6 \gl^2 (1 + 3 \gl^2)}
\le \frac2\gl \Bigpar{1 - \frac\gl2} = \frac2\gl - 1.
$$
Hence, \eqref{g'} shows that $g_{12}$ is decreasing on $[1/3,1]$, and thus
$$
g_{12}(\gl) \le g_{12}(0.42) < 1,
\qquad 0.42 \le \gl \le1.
$$
Finally,
$$
g_{12}(\gl) \le \frac{1}{6} \frac{e^{\gl}}{\gl^2}\leq \frac{e}{6} < 1,
\qquad 1 \le \gl \le 2.
$$
Combining these estimates we find that if $K \ge \max(12, 2L^{-2}e^L)$,
then $f_K(\gl) \le g_K(\gl) < 1$ whenever $0.42 \le \gl \le L$, while
$f_K(\gl) \le f_1(\gl) \le 1$ for $0 \leq \gl \le 0.42$, and thus
\eqref{fkl} holds when $\lambda \geq 0$.

We have also shown that $K = 12$ will do for $L \le 2$ and $\gl\ge0$;
since $2 L^{-2} e^L$ is increasing for $L \ge 2$, and thus less than $12$
for $2 \le L < L_0$ but larger than $12$ for $L > L_0$,
\refT{T:KL} for $\gl\ge0$ follows.

For $\lambda \leq 0$, we again use Taylor's formula for small $|\gl|$;
arguing as above we see that \eqref{fkl} holds for $\gl\le0$ provided
\begin{equation}
\label{thomas}
L\sup_{-L \le\gl\le 0}\bigpar{- f_K'''(\gl)}
\le2K-\gss.
\end{equation}
It is easily checked numerically that
$\max_{0\le u\le 1} u(1-u)g(u)<0.033$.
It follows that
$$
X(u,\gl)
\geq (- \eta^3 - 0.396 K - 3 K^2 L) e^{\eta L},
\qquad -L \leq \lambda \leq 0.
$$
Hence, \eqref{thomas} holds and
\eqref{fkl} is satisfied
for $\lambda \leq 0$ provided
$$
(\eta^3 + 0.396 K + 3 K^2 L) L e^{\eta L} \leq 2 K - \sigma^2.
$$
It is readily checked that this holds for $K = 0.5$ and $L = 0.62$.

For larger $L$ we argue as follows.
The function $h(u):=g(u)+4\eta u(1-u)$ satisfies
$$
h''(u)=\frac2{u(1-u)} -8\eta \ge 8-8\eta>0,
\qquad 0<u<1.
$$
Hence $h$ is convex, and since $h'(1/2)=0$,
$$
h(u)\ge h(\tfrac12)=0,
\qquad 0\le u\le1.
$$
Consequently, if $\gl\le0$ and $K|\gl|\ge2\eta$, then
$$
\gl g(U)-2K\gl^2 U(1-U) \le \gl h(U) \le0
$$
and thus $f_K(\gl)\le1$.
Choosing $K=2\eta/0.62<1.247$, this shows $f_K(\gl)\le1$ for
$\gl\le-0.62$, while
$f_K(\gl)\le f_{0.5}(\gl)\le 1$ for $-0.62\le \gl \le0$
by the preceding case.

This completes the proof of both theorems.
\end{proof}

If we just want a bound  on $\psi_Y$, \eqref{psi1} and Theorems
\refand{T:KL}{T:KL_neg}
can be stated more simply as follows (ignoring the better bounds
obtained for small $\gl$).
\begin{corollary}
\label{C:Ymgfbds}
With $L_0$ as
in \refT{T:KL},
\begin{equation*}
\psi_Y(\gl) \le
\begin{cases}
e^{1.25 \gl^2},            & \gl \le 0,\\
e^{12 \gl^2},           & 0\le \gl \le L_0,\\
e^{2 e^{\gl}}, & \gl \ge L_0.
\end{cases}
\end{equation*}
In particular, $\psi_Y(\gl)\le \exp(\max(12\gl^2,2e^\gl))$.
\nopf
\end{corollary}

The bound $e^{2e^\gl}$ is very large even for moderately large $\gl$,
but the next result shows that $\psi_Y(\gl)$ really is of essentially
this size. In particular, it follows that $\ln\ln\psi_Y(\gl)\sim\gl$
as $\gl\to+\infty$.
\begin{theorem}
\label{T:mgf_lower}
If $\gamma<2/e$, then for sufficiently large $\gl$,
$$
\psi_Y(\gl)\ge \exp(\gamma\gl^{-1}e^\gl).
$$
\end{theorem}
\begin{proof}
Since a \mgf{} is convex and $\psi_Y'(0)=\EE Y=0$, $\psi_Y$ is
increasing on $[0,\infty)$.
Moreover, $g$ is decreasing on $[0,1/2]$.
Hence, if $0\le\gd\le1/2$,
the integrand in \eqref{psiinteq} with $Z=Y$ is for $0\le u\le \gd$
at least $\psi_Y(0)\psi_Y((1-\gd)\gl)e^{\gl g(\gd)}$
and the same holds for $1-\gd\le u\le 1$ by symmetry.
Consequently,
\begin{equation}
\label{magnus}
\psi_Y(\gl)
\ge2\int_0^\gd
 \psi_Y(u \lambda)\,\psi_Y((1 - u) \lambda)\,e^{ \lambda g(u)}\,du
\ge 2\gd \psi_Y((1-\gd)\gl)e^{\gl g(\gd)},
\qquad 0\le\gd\le1/2.
\end{equation}
Let $a>1/2$ be a constant to be determined later and choose
$\gd := a e^{-\gl}$. Then $g(\gd)=1-O(\gl e^{-\gl})$ and thus by
\eqref{magnus}, for
$\gl\ge\ln(2a)$,
$$
\psi_Y(\gl)\ge 2a e^{-O(\gl^2 e^{-\gl})} \psi_Y(\gl-a\gl e^{-\gl}).
$$
If $0<\eps<2a$, there thus exists $A$ such that for
$\gl\ge A$,
$$
\psi_Y(\gl)\ge (2a-\eps)\psi_Y(\gl-a\gl e^{-\gl}).
$$

Given $\gl \geq A$, let
$\gl_0 := \gl$ and define inductively
$\gl_{n+1} := \gl_n-a\gl_ne^{-\gl_ n}$, $n\ge0$.
Let $N$ be the smallest integer with $\gl_N<A$.
Then
$
\psi_Y(\gl_{n})\ge (2a-\eps)\psi_Y(\gl_{n+1})
$,
$n=0,\dots,N-1$, and thus
$$
\psi_Y(\gl)=\psi_Y(\gl_{0})
\ge (2a-\eps)^N\psi_Y(\gl_N)
\ge (2a-\eps)^N.
$$
It remains to estimate $N$ from below. Since $e^x$ is increasing,
$$
\int_{\gl_{n+1}}^{\gl_n} e^x\,dx
\le e^{\gl_n}(\gl_n-\gl_{n+1})
=a\gl_n \le a\gl
$$
and thus
$$
Na\gl
\ge \int_{\gl_N}^{\gl_0} e^x\,dx
\ge \int_{A}^{\gl} e^x\,dx
=e^\gl-e^A.
$$
Consequently,
$$
\ln\psi_Y(\gl) \ge N \ln(2a-\eps)
\ge \frac{\ln(2a-\eps)}{a} \gl^{-1}(e^\gl-e^A),
\qquad \gl \ge A.
$$
We choose $a=e/2$, which maximizes $\ln(2a)/a$. Then $\ln(2a)/a=2/e$,
we may choose $\eps$ so small that $\ln(2a-\eps)/a>\gamma$, and
the result follows.
\end{proof}

As is well known, bounds on the \mgf{} yield bounds on the tails of
the distribution.

\begin{theorem}\label{T:tail}
If $y\ge 2e^{L_0}=12 L_0^2 \doteq 302.1$, then 
$$
\PP(Y\ge y) \le \exp\bigpar{-y(\ln y-1-\ln 2)}.
$$
\end{theorem}

\begin{proof}
For $\gl\ge L_0$, by \refC{C:Ymgfbds},
$$
\PP(Y\ge y) \le e^{-\gl y}\EE e^{\gl Y} \le \exp(2e^\gl-y\gl),
$$
and the result follows by taking $\gl=\ln(y/2)$.
\end{proof}
\begin{remark}
The same estimate holds for every $Z_n$ provided, say,
$\psi_{Z_0}(\gl)\le\exp(12\gl^2)$; for example, when
$Z_0\equiv0$.
\end{remark}

\refT{T:mgf_lower} suggests
that the true size of $\PP(Y\ge y)$ is (for large~$y$) not
much smaller than the upper bound in \refT{T:tail}.
Indeed,
Knessl and Szpankowski \cite{KnSz}
have found (assuming an as yet unverified
regularity hypothesis) a much more precise formula for the asymptotics
of $\PP(Y\ge y)$ which is of the order
$\exp\bigpar{-y [\ln y + \ln\ln y + O(1)]}$.

For the left tail,
\refC{C:Ymgfbds} similarly implies
$\PP(Y\le y)\le \exp(-y^2/5)$ for $y\le0$, but this result is much
weaker than the doubly exponential decay found by
Knessl and Szpankowski \cite{KnSz}.

\section{Geometric rate of convergence for moment generating functions}
\label{S:mgf2}

\begin{theorem}\label{T:mgf}
Suppose that $Z_0$ has mean zero and an everywhere finite \mgf{}~$\psi_{Z_0}$.
Then $\psiznl\to\psiyl$ at a geometric rate for every
fixed $\gl \in \RR$.
Explicitly, if $L \geq 0$ and $K_L$ are such that
\eqref{psi1} and \eqref{psi2} hold, and if moreover
\begin{equation}\label{psi3}
\psi_{Z_0}(\gl)\le e^{K_L\gl^2},
\qquad |\gl|\le L,
\end{equation}
then, for every $n\ge0$ and $|\gl|\le L/2$,
\begin{equation*}
\begin{split}
|\psiznl-\psiyl|
&\le (\Var\,Z_0 + \gss)^{1/2}|\gl|
 \bigpar{\psizn(2\gl)+\psiy(2\gl)}^{1/2} (2/3)^{n/2}\\
&\le 2^{1/2}(\Var\,Z_0 + \gss)^{1/2}|\gl|
  e^{2K_L\gl^2}  (2/3)^{n/2}.
\end{split}
\end{equation*}
\end{theorem}
Of course, if the hypotheses in the first sentence of the theorem's statement
are met, then, given $L \geq 0$, \eqref{psi3} holds for some $K_L <
\infty$, which by \refT{T:KL} may be chosen so large that \eqref{psi1}
and \eqref{psi2} hold, too.

\begin{proof}
By \eqref{psi2} and induction, the estimate \eqref{psi3} holds for
every $\psizn$.
Fix $n \geq 0$ and consider the optimal $d_2$-coupling of (the laws
of)~$Z_n$ and~$Y$.
Then for $\lambda \in [- L / 2,\,L / 2]$ we have,
using the mean value theorem and the Cauchy--Schwarz inequality,
\begin{equation*}
\begin{split}
\left| \EE\,e^{\lambda Z_n} - \EE\,e^{\lambda Y} \right|
&\leq \EE \left| e^{\lambda Z_n} -e^{\lambda Y} \right|\\
&\leq \EE \left(|\lambda| |Z_n - Y|  e^{\max(\lambda Z_n,\lambda
Y)}\right)\\
&\leq |\lambda| \bigpar{\EE|Z_n - Y|^2}^{1/2}
\Bigpar{ \EE\,e^{2\max(\lambda Z_n,\lambda Y)}}^{1/2} \\
&\leq |\lambda| \bigpar{\EE|Z_n - Y|^2}^{1/2}
\Bigpar{ \EE\,e^{2\lambda Z_n}+\EE\,e^{2\lambda Y}}^{1/2}.
\end{split}
\end{equation*}
By the optimality of the coupling and
\eqref{d2conv},
$$
\bigpar{\EE|Z_n - Y|^2}^{1/2}  =d_2(F_n,F)
\le (\Var\,Z_0 + \gss)^{1/2}(2/3)^{n/2},
$$
and by~\eqref{psi1} and~\eqref{psi3} for~$\psi_{Z_n}$
$$
\EE\,e^{2\lambda Z_n} + \EE\,e^{2\lambda Y} \le 2 e^{K_L(2\gl)^2},
$$
whence the result follows.
\end{proof}

Note further that the operator $\psi_Z\mapsto\psi_{S Z}$
given by~\eqref{psiinteq} is monotone, in the specific sense that if
$\psi_Z(\gl)\le\psi_W(\gl)$ for $|\gl|\le L$,
then also $\psi_{S Z}(\gl)\le\psi_{S W}(\gl)$ for $|\gl|\le L$.
In particular, by induction,
if $\psi_{Z_0}(\gl)\le\psi_{Z_1}(\gl)$ for $|\gl|\le L$,
then  $\psi_{Z_n}(\gl)$ increases monotonically to its limit
$\psiyl$ for $|\gl|\le L$.  Likewise,
 if $\psi_{Z_0}(\gl)\ge\psi_{Z_1}(\gl)$ for $|\gl|\le L$,
then  $\psi_{Z_n}(\gl)$ decreases to $\psiyl$ for $|\gl| \leq L$.

We give two simple special cases.

\begin{corollary}
Suppose that $Z_0\equiv0$.
Then $\psiznl$ increases monotonically to $\psiyl$ for every
fixed $\gl$.
If $L \geq 0$ and $K_L$ are such that
\eqref{psi1} and \eqref{psi2} hold,
then, for every $n\ge0$ and $|\gl|\le L/2$,
\begin{equation*}
\begin{split}
0\le\psiyl-\psiznl
&\le 2^{1/2} \gs |\gl| \left( \psiy(2\gl) \right)^{1/2} (2/3)^{n/2}\\
&\le 2^{1/2}\gs|\gl|  e^{2K_L\gl^2}  (2/3)^{n/2}.
\end{split}
\end{equation*}
\end{corollary}
\begin{proof}
Since $\EE\,Z_1=0$, by Jensen's inequality
$\psi_{Z_1}(\gl) \ge 1= \psi_{Z_0}(\gl)$, and the monotonicity
follows.
In particular, $\psizn(2\gl)\le\psiy(2\gl)$, and
since \eqref{psi3} trivially is satisfied,
the result follows
from \refT{T:mgf}.
\end{proof}

\begin{corollary}
Suppose  $L \geq 0$ and $K_L$ are such that
\eqref{psi1} and \eqref{psi2} hold,
and let $Z_0 \sim N(0,2K_L)$.
Then $\psiznl$ decreases monotonically to $\psiyl$ for every
fixed $\gl$ with $|\gl|\le L$,
and, for every $n\ge0$ and $|\gl|\le L/2$,
\begin{equation*}
\begin{split}
0\le\psiznl-\psiyl
&\le (4 K_L + 2 \gss)^{1/2} |\gl| \left( \psizn(2\gl) \right)^{1/2}
(2/3)^{n/2}\\
&\le (4K_L+2\gss)^{1/2}|\gl|  e^{2K_L\gl^2}  (2/3)^{n/2}.
\end{split}
\end{equation*}
\end{corollary}
\begin{proof}
Since $\psi_{Z_0}(\gl)=e^{K_L\gl^2}$,
\eqref{psi2} yields
$\psi_{Z_1}(\gl) \le \psi_{Z_0}(\gl)$, and the monotonicity
follows.
The estimate thus follows from \refT{T:mgf}.
\end{proof}

\section{On numerical calculations}\label{S:numerics}

The preceding results make it possible, in principle at least, to
calculate the density, distribution, characteristic, and moment
generating functions of $Y$ numerically, with
provable arbitrarily high accuracy.

To begin, the results of earlier sections show that it suffices to start with a
suitable $\Lc(Z_0)$, for example unit mass at~$0$ or a normal distribution,
and then calculate the corresponding quantity for $Z_n$, for a large~$n$
that can be determined.
The distribution of $Z_n$ can be calculated recusively; for the
characteristic and moment generating functions we have the recurrence
relations~\eqref{phiinteq} and~\eqref{psiinteq}, while for the density
functions we have the following recursion:
\begt\label{T:inteq}
If $n \geq 0$ is arbitrary and~$Z_0$ has a bounded continuous density
function~$f_0$,
or if~$Z_0$ is arbitrary and $n \ge 3$, then~$Z_n$ and~$Z_{n + 1}$
have bounded continuous density functions~$f_n$ and~$f_{n + 1}$ satisfying
the identity
\begin{equation}\label{densrec}
f_{n+1}(x) = \int^1_{u = 0}\!\int_{z \in \RR}\!f_n(z)\,f_n
\left( \frac{x - g(u) - (1 - u) z}{u} \right)
\frac{1}{u}\,dz\,du,\qquad x \in \RR,
\end{equation}
with~$g(\cdot)$ given by~\eqref{gu}.
\ent
\begin{proof}
Our proof
(similar to that of Theorem~4.1 in~\cite{SJFill1})
is by induction on $n \geq 0$ in the first case and on $n \geq 3$, using
\refT{T:density} to get started, in the second case.  We may therefore assume
as our induction hypothesis that~$f_n$ is bounded and continuous.
It is easily checked that, for each $0 < u < 1$, the inner integral
\begin{equation*}
h_u(x):= \int_{z \in \RR}\!f_n(z)\,f_n \left( \frac{x - g(u) - (1 - u)
z}{u} \right)
\frac{1}{u}\,dz
\end{equation*}
is a density function for
the random variable
\begin{equation}
\label{fixurv}
u Z_n + (1 - u) Z_n^* + g(u),
\end{equation}
and, using dominated convergence, that~$h_u$ is bounded and
continuous.
Indeed,  $h_u(x) \le (\sup f_n) / u$, and since
$h_u=h_{1-u}$ by symmetry in~\eqref{fixurv},
$h_u(x) \le 2\sup f_n$, uniformly in $u$ and $x$.
It follows, by dominated convergence again, that
$x \mapsto f_{n + 1}(x) = \int_0^1\!h_u(x)\,du$
is a bounded continuous density for~$Z_{n+1}$.
\end{proof}

The integrals in~\eqref{phiinteq}, \eqref{psiinteq},
or~\eqref{densrec} have to be
computed numerically---as does the integral of~$f_n$ to get~$F_n$---but
that can be done with arbitrary precision since the results above
provide bounds for the
integrands and their derivatives. [The function $g(u)$ has an
unbounded derivative as $u\to0$ or $u\to1$, but that can be handled by
truncating the interval.] Consequently, to calculate $\phizn(t)$ with
given precision for a given $t$, it suffices to know $\phi_{Z_{n-1}}(t_k)$
with another given precision for a finite number of points $t_k$,
which can be done recursively.
(However, a brute force recursion along these lines seems to require
too many numerical integrations to be practical if we want reasonably
high provable accuracy.)

\begin{remark}
To calculate the density $f_n$ numerically, it might be better to
compute $\phi_{Z_n}$ recursively by \eqref{phiinteq} and then use
\eqref{Fourierinv}, instead of using
the recursion \eqref{densrec} directly.  This is both because
\eqref{densrec} is a double integral and because we have the simple
bounds
$|\phi_n|\le1$ and $|\phi_n^{(k)}|\le \EE|Z_n|^k$, $k\ge1$.
\end{remark}

\section{The metrics $d_p$ and a lower bound on $d_2(F_n,
F)$}\label{S:lower}

At~\eqref{d2conv} we recalled \Roesler's fundamental result
$$
d_2(F_n, F) = O(\rho^n)
$$
with $\rho := (2 / 3)^{1 / 2}$.  The question naturally arises as to whether
there is a lower bound that matches at least to the extent that
$$
d_2(F_n, F) = \Omega(r^n)
$$
for some $r > 0$.  Of course, the answer is negative without any
further restrictions,
since if $F_0 = F$ then $F_n = F$ for every~$n$.  However, our main
result of this
section asserts that this is the only exception, at least among
distributions~$F_0$
with finite moments of all orders:

\begt\label{T:LB}
If $F_0 \neq F$ has finite moments of all orders, then there exists $r
> 0$ (depending
on~$F_0$) so that
$$
d_2(F_n, F) = \Omega(r^n).
$$
\ent

Our arguments for \refT{T:LB} will require use of metrics~$d_p$
generalizing~\eqref{d2}.  So we will warm up in \refS{S:dp} by
recalling the definition
of, and two useful facts about, $d_p$ and in \refS{S:UBp} by extending
the upper bound
result~\eqref{d2conv} to~$d_p$ for $p \geq 1$.  Then in \refS{S:LB} we
will prove a
sharpened version of \refT{T:LB} (namely, \refT{T:sharper}).

\subsection{The metrics~$d_p$} \label{S:dp}

For real $1 \leq p < \infty$, let $\| X \|_p := \left( \EE|X|^p
\right)^{1 / p}$ denote
the $L^p$-norm, and let~$d_p$ denote the metric on the space of probability
distributions with finite $L^p$-norm defined by
$$
d_p(F, G) := \min \| X - Y \|_p,
$$
taking the minimum, as at~\eqref{d2}, over all couplings of $\Lc(X) =
F$ and $\Lc(Y) =
G$.  It is worth noting
that there is a coupling [namely, $X =
F^{-1}(U)$ and $Y = G^{-1}(U)$ for~$U$ uniform and a suitable
definition of the inverse
probability transform~$F^{-1}$] that achieves the minimum
simultaneously for each $1
\leq p < \infty$ (assuming~$F$ and~$G$ have finite moments of all
orders):\ see~\cite{Cambanis}. 

We begin with two elementary facts that will be useful later.
The proof of the first fact (\refL{L:decrease}) shows that~$S$
is a contraction for the $d_p$-metric.
\begl\label{L:decrease}
Consider real $1 \leq p < \infty$.  The $d_p$-distance from the
limiting {\tt Quicksort} distribution~$F$ does not increase when the
operator~$S$ of~\eqref{operator}
is applied.  Therefore, $d_p(F_n, F)$ is nonincreasing, and hence
bounded, in~$n$ if
$\EE\,|Z_0|^p < \infty$.
\enl
\begin{proof}
With a slight abuse of notation, we find, for~$Z$ with any law,
\begin{equation}
\label{dpineq}
d_p(SZ, Y) = d_p(SZ, SY) \leq \| U (Z - Y) + (1 - U) (Z^* - Y^*) \|_p,
\end{equation}
coupling $(Y, Z)$ optimally and $(Y^*, Z^*)$ optimally and choosing
$U$, $(Y, Z)$, and $(Y^*, Z^*)$ to be independent.  In calculating the
$L^p$-norm value on the right in~\eqref{dpineq}, condition on~$U$ and then
use subadditivity of $L^p$-norm together with independence to bound that value
by $\| Z - Y \|_p = d_p(Z, Y)$.  This establishes the first
assertion:\ $d_p(SZ, Y) \leq
d_p(Z, Y)$.

Therefore, $d_p(F_n, F) = d_p(S^n Z_0, Y)$ is nonincreasing, and hence
bounded by $d_p(Z_0, Y)$, which is bounded by $\| Z_0 \|_p + \| Y \|_p
< \infty$
if $\EE\,|Z_0|^p < \infty$.
\end{proof}

\begin{remark}\label{R:dpinfty}
Conversely, if $\|SZ\|_p<\infty$, then $\EE|uZ+(1-u)Z^*|^p<\infty$ for
some $u\in(0,1)$, and thus $\EE|Z|^p<\infty$ too.
Hence, if $\EE|Z_0|^p=\infty$, then $\EE|Z_n|^p=\infty$ and
$d_p(Z_n,Y)=\infty$ for all $n$.
\end{remark}

\begl\label{L:inter}
For real $2 < p < q < \infty$ we have, for any~$F$ and~$G$,
$$
d_p(F, G) \leq
d_2^{\frac{2 (q - p)}{p (q - 2)}}(F, G)
  \times d_q^{\frac{q (p - 2)}{p (q - 2)}}(F, G).
$$
\enl

\begin{proof}
Using
the common optimal coupling for $d_2$ and $d_q$,
this is immediate from the inequality
$$
\| X \|_p \leq \| X \|_2^{\frac{2 (q - p)}{p (q - 2)}} \| X
\|_q^{\frac{q (p - 2)}{p (q - 2)}},
$$
which in turn follows from the
fact~\cite[Exercise~4(b) of Chapter~3]{RudinRealComplex}
that $\ln \| X \|^p_p$ is convex in $p \in (0, \infty)$.
\end{proof}

\subsection{Geometric rate of convergence in each metric~$d_p$}\label{S:UBp}

Under suitable conditions, we can establish a geometric rate of convergence
for $d_p(F_n, F)$ for any real $1 \leq p < \infty$.
We begin with an elementary lemma.

\begin{lemma} \label{L:xp}
If $p\ge2$, then for all $x,y\ge0$,
$$
(x+y)^p\le x^p +y^p +c_p(x^{p-1}y+xy^{p-1}),
$$
where $c_p := p(p-1)2^{p-2}$.
\end{lemma}
\begin{proof}
\begin{equation*}
\begin{split}
(x+y)^p-x^p-y^p
&=\int_0^yp\bigpar{(x+t)^{p-1}-t^{p-1}}\,dt
=\int_{t=0}^y\int_{u=0}^xp(p-1)(t+u)^{p-2}\,du\,dt\\
&\le p(p-1)xy(x+y)^{p-2}
\le p(p-1) xy 2^{p-2}(x^{p-2}+y^{p-2}).
\end{split}
\qedtag
\end{equation*}
\noqed
\end{proof}

\begt\label{T:UBp}
Let $p_0\doteq 6.557$ be the largest positive solution to
\begin{equation}\label{p0}
\parfrac{2}{p_0+1}^{1/p_0} = \parfrac 23^{1/2}
\end{equation}
and
let, for any $\eps>0$,
\begin{equation}\label{betap}
\beta_p:=\begin{cases}
\bigpar{\frac 23}^{1/2}, & 1\le p < p_0,\\
\bigpar{\frac 23}^{1/2}+\eps, & p = p_0,\\
\bigpar{\frac{2}{p+1}}^{1/p} & p>p_0.
\end{cases}
\end{equation}
Thus, for $p\ge 2$ except $p=p_0$,
$\beta_p=\max\bigpar{(2/3)^{1/2},(2/(p+1))^{1/p}}$.
Then, for any~$Z_0$ with zero mean and finite variance,
and every $p\ge1$ such that $\EE|Z_0|^p<\infty$,
there exists a constant $\alpha_p < \infty$ [depending on $\Lc(Z_0)$]
such that
\begin{equation}
\label{UBp}
d_p(Z_n, Y) \leq \alpha_p \beta^n_p.
\end{equation}
\ent

\begin{proof}
First we note that \eqref{p0} can be written $(3/2)^{p_0/2}=(p_0+1)/2$.
One root of this equation is $2$, and since $(3/2)^{p/2}$ is convex,
with derivative less than $1/2$ at $p=2$, it follows that the equation
has two positive roots,
$2$ and $p_0>2$, and that $(2/3)^{p/2} >2/(p+1)$ for $2<p<p_0$, while
$(2/3)^{p/2} <2/(p+1)$ for $p>p_0$.  

Next
we note that~\eqref{UBp} holds for $p\le2$, with
$$
\alpha_p := (\Var\,Z_0 + \sigma^2)^{1 / 2}, \qquad p\le2,
$$
by~\eqref{d2conv} and the inequality $d_p \leq d_2$, $p\le2$.
We then proceed by induction on~$\floor p$.  For the induction step,
suppose that $p >2$
and that~$Z_0$ has zero mean
and satisfies $\| Z_0 \|_p < \infty$.
By the induction hypothesis, there exist constants $0 < \alpha_q <
\infty$, $1\le q\le  p - 1$, such that
\begin{equation}
\label{IH}
d_q(Z_n, Y) \leq
 \alpha_q \beta^n_q\text{\ \ for all $q\le p- 1$ and $n \geq 0$.}
\end{equation}
Using our usual coupling of~$Z_n$ and~$Y$ in terms of the optimal
coupling of~$Z_{n - 1}$ and~$Y$, we find easily
by \refL{L:xp},
for $n \geq 1$
[with $(Z_{n-1},Y)$, $(Z_{n-1}^*,Y^*)$, and $U$ independent],
\begin{equation*}
\begin{split}
d^p_p(Z_n, Y)
  &\leq \EE\bigl|U(Z_{n-1}-Y)+(1-U)(Z_{n-1}^*-Y^*)\bigr|^p \\
  &\leq \EE\bigl(U|Z_{n-1}-Y|+(1-U)|Z_{n-1}^*-Y^*|\bigr)^p \\
  &\leq \EE \bigpar{U^p|Z_{n-1}-Y|^p}+\EE\bigpar{(1-U)^p|Z_{n-1}^*-Y^*|^p} \\
& \qquad  +c_p\EE \bigpar{U^{p-1}(1-U) |Z_{n-1}-Y|^{p-1}|Z_{n-1}^*-Y^*|} \\
& \qquad   +c_p\EE \bigpar{U(1-U)^{p-1}  |Z_{n-1}-Y||Z_{n-1}^*-Y^*|^{p-1} }\\
&=\frac2{p+1} d^p_p(Z_{n - 1}, Y)
+\frac{2c_p}{p(p+1)} d_1(Z_{n - 1}, Y)\,d^{p -1}_{p - 1}(Z_{n - 1}, Y).
\end{split}
\end{equation*}
So by induction on~$n$ it follows that
\begin{equation*}
d_p^p(Z_n,Y)\le
\parfrac{2}{p + 1}^n d_p^p(Z_0,Y) + \frac{c_p}p\sum_{i = 0}^{n-1}
 \parfrac{2}{p + 1}^{n - i} d_1(Z_i,Y)d_{p-1}^{p-1}(Z_i,Y).
\end{equation*}
By the induction hypothesis \eqref{IH} this yields, for some
$a_1,a_2<\infty$ (depending  on~$p$),
\begin{equation}
\label{sofie}
d_p^p(Z_n,Y)\le
a_1 \parfrac{2}{p + 1}^n  + a_2\sum_{i = 0}^{n-1}
 \parfrac{2}{p + 1}^{n - i}\bigpar{\beta_1\beta_{p-1}^{p-1}}^i.
\end{equation}

Let  $\gamma := \beta_1 \beta_{p-1}^{p-1}$.
We break our treatment into three cases:
\begin{enumerate}
\item
If $\gamma>2/(p+1)$, we write the sum in \eqref{sofie} as
$$
\gamma^n\sum_{i=0}^{n-1}\left(\frac2{p+1}\gamma^{-1}\right)^{n-i}
<\Bigpar{1-\frac2{p+1}\gamma^{-1}}^{-1} \gamma^n,
$$
and thus \eqref{sofie} shows that \eqref{UBp} holds  with $\beta_p^p=\gamma$.
\item
If $\gamma<2/(p+1)$, we write the sum in \eqref{sofie} as
$$
\parfrac2{p+1}^n\sum_{i=0}^{n-1}\left(\gamma\frac{p+1}2\right)^{i}
<\Bigpar{1-\gamma\frac{p+1}2}^{-1} \parfrac2{p+1}^n,
$$
and thus \eqref{UBp} holds  with $\beta_p^p=2/(p+1)$.
\item
If $\gamma=2/(p+1)$,  the sum in \eqref{sofie} equals  $n\bigpar{2/(p+1)}^n$.
Consequently, \eqref{UBp} holds  with any $\beta_p>\bigpar{2/(p+1)}^{1/p}$.
\end{enumerate}

It remains to verify that this yields the $\beta_p$ given in \eqref{betap}.

First, if $2<p <p_0$, then
the induction hypothesis yields
$\gamma=\beta_1\beta_{p-1}^{p-1} = (2/3)^{p/2} > 2/(p+1)$,
so case (i) gives
$\beta_p^p=(2/3)^{p/2}$.
Similarly, for $p=p_0$,
$\gamma= (2/3)^{p/2}= 2/(p+1)$ and (iii) shows that any
$\beta_p>\bigpar{2/(p+1)}^{1/p} = (2 / 3)^{1 / 2}$ 
will do.

For $p_0<p <p_0+1$, we have
$\gamma=\beta_1\beta_{p-1}^{p-1} = (2/3)^{p/2} < 2/(p+1)$,
so case (ii) yields
$\beta_p^p=2/(p+1)$.
The same applies for $p=p_0+1$, since
again $(2/3)^{p/2} < 2/(p+1)$ and 
we thus may choose $\eps$ so small that
$\gamma=(\frac23)^{1/2}\bigpar{(\frac23)^{1/2}+\eps}^{p-1}<2/(p+1)$.

Finally, for $p>p_0+1$, we have
$$
\gamma=\beta_1\beta_{p-1}^{p-1} = \biggpar{\frac 23}^{\!1/2}\,\frac2p
< \frac2 {p+1},
$$
since $p/(p+1)$ is increasing and equals $(2/3)^{1/2}$ when
$p=(\sqrt{3/2}-1)^{-1} = 2(\sqrt6-2)^{-1}=\sqrt6+2<5<p_0$.
Hence case (ii) applies.
\end{proof}

\subsection{Lower bounds}\label{S:LB}
The main goal of this subsection is to establish \refT{T:LB}, or
rather the sharper
\refT{T:sharper} below.  Since (as noted in \refS{S:mgf1}) the limiting
{\tt Quicksort}
distribution~$F$ has everywhere finite \mgf, it is uniquely determined
by its moments.
Hence if $F_0 \neq F$ has finite moments of all orders, then
$\EE\,Z^j_0 \neq \EE\,Y^j$
for some integer $j \geq 1$.

\begt\label{T:sharper}
Suppose $F_0 \neq F$ has finite moments of all orders, and let~$p$
be the smallest positive integer such that
$\EE\,Z^p_0 \neq \EE\,Y^p$.
Then, for any $0 < r < \bigpar{\sfrac{2}{p + 1}}^{p / 2}$,
$$
d_2(F_n, F) = \Omega(r^n).
$$
(The implicit multiplicative constant depends on both~$F_0$ and~$r$.)
\ent

\begin{remark}\label{R:d12}
The
cases $p=1$ and $p=2$ are a bit special, and in these cases
we claim that
\refT{T:sharper} holds even
with $r=\bigpar{\sfrac{2}{p + 1}}^{p / 2}$,
\ie{}, with $r = 1$ and $r = 2/3$, respectively, and without the assumption
that $F_0$ has finite moments.

We may
and shall assume that $\EE\,Z_0^2 < \infty$,
since otherwise $d_2(F_n, F) = \infty$. 

First, 
$p=1$ when 
$\EE\,Z_0\neq \EE\,Y=0$; in this case 
$Z_n$
converges in distribution to $Y+\EE\,Z_0$ and not to $Y$, and thus
$\inf d_2(Z_n,Y)>0$,
\ie, the theorem holds with $r=1$.
Indeed, 
we have the sharper result that
$$
d_2(Z_n, Y)= d_2(Y + \EE\,Z_0, Y) + O\bigpar{d_2(Z_n, Y + \EE\,Z_0)}
=|\EE\,Z_0|+O\bigpar{(2/3)^{n/2}}.
$$

Next, $p=2$ when $\EE\,Z_0=0$ but $\Var\,Y \neq \Var\,Z_0$;
in this case \refT{T:LBp} shows that the result holds with
$r=2/3$. Even in this case we have a gap between the lower
bound $\Omega\bigpar{(2/3)^n}$ and \Roesler's upper bound
$O\bigpar{(2/3)^{n/2}}$; 
it is an open problem to find the
rate of approximation more precisely.
\end{remark}
We prove
\refT{T:sharper} using the following \refT{T:LBp},
which is a similar lower bound for the $d_p$-metric.

\begt\label{T:LBp}
Let
$p \geq 1$ be an integer, and suppose that $\EE\,Z^j_0 = \EE\,Y^j$ for
integers $1
\leq j \leq p - 1$, and that $\EE\,Z^p_0$ exists and is finite but
fails to equal
$\EE\,Y^p$.  Then
$$
d_p(F_n, F) = \Omega \left( \bigpar{\sfrac{2}{p + 1}}^n \right).
$$
\ent

\refT{T:LBp} is, in turn, a simple consequence of the following two
elementary lemmas.

The first lemma demonstrates a sense in which the value of~$p$
in Theorems \refand{T:sharper}{T:LBp}
persists from~$F_0$ to each~$F_n$; the second gives a
general lower bound
on~$d_p$ in terms of discrepancy in $p$th moments.

\begl\label{L:persist}
Let $p \geq 1$ be an integer, and suppose for $n = 0$ that $\EE\,Z^j_n
= \EE\,Y^j$
for integers $1 \leq j \leq p - 1$, and that $\EE\,Z^p_n$ exists and
is finite but fails
to equal $\EE\,Y^p$.  Then for every $n \geq 0$ the same is true and, moreover,
$$
\EE\,Z^p_n - \EE\,Y^p
= \bigpar{\sfrac{2}{p + 1}}^n \left( \EE\,Z^p_0 - \EE Y^p \right).
$$
\enl

\begin{proof}
If $\EE|Z|^m<\infty$,
then, with $Z$, $Z^* \Leq Z$, and  $U$ independent, by
\eqref{operator} and a trinomial expansion we have
\begin{equation*}
\begin{split}
\EE(SZ)^m
&=\sum_{j+k\le m} \frac{m!}{j!\,k!\,(m-j-k)!}
\EE\bigpar{U^j Z^j(1-U)^k (Z^*)^k g(U)^{m-j-k}} \\
&=\sum_{j+k\le m} \frac{m!}{j!\,k!\,(m-j-k)!}
\EE\bigpar{U^j(1-U)^kg(U)^{m-j-k}}\,\EE Z^j \,\EE Z^k .
\end{split}
\end{equation*}

We apply this with $m=1,\dots,p$ for both $Z=Z_{n-1}$ and $Z=Y$, and
note that by induction on~$n$, all terms in the sum with $j\le p-1$
and $k\le p-1$ coincide for the two choices of $Z$. Hence,
$\EE Z_n^m=\EE Y^m$ for $1\le m\le p-1$ and
$$
\EE Z_n^p -\EE Y^p
=  \bigpar{\EE U^p+\EE (1-U)^p}\bigpar{\EE Z_{n-1}^p -\EE Y^p }
=  \frac2{p+1}\bigpar{\EE Z_{n-1}^p -\EE Y^p },
$$
and the result follows.
\end{proof}

\begl\label{dpmoments}
Let $p \geq 1$ be an integer.  Then, for any~$F$ and~$G$,
$$
d_p(F, G) \geq
\frac{| \EE\,X^p - \EE\,Y^p |}
{\sum_{j = 0}^{p - 1} \| X \|^{p - 1 -j}_p\,\| Y \|^j_p}
$$
with $X \sim F$ and $Y \sim G$ (and $0^0 := 1$).
\enl

\begin{proof}
Let $(X, Y)$ be an optimal coupling of~$F$ and~$G$.  If $p = 1$, then
$$
d_1(F, G) = \| X - Y \|_1 = \EE\,|X - Y| \geq | \EE\,X - \EE\,Y |,
$$
as desired.  If $p \geq 2$, we employ the factorization
$$
X^p - Y^p = (X - Y) \sum_{j = 0}^{p - 1} X^{p - 1 - j} Y^j,
$$
whence
\begin{eqnarray}
|\EE\,X^p - \EE\,Y^p|
  &\leq& \EE |X^p - Y^p| \nonumber \\
  &\leq& \| X - Y \|_p\,\left\|
  \sum_{j = 0}^{p - 1} X^{p - 1 - j} Y^j \right\|_{p / (p -
           1)} \nonumber \\
\label{intermediate}
  &\leq& d_p(F, G) \sum_{j = 0}^{p - 1} \left\| X^{p - 1 - j} Y^j
  \right\|_{p / (p - 1)},
\end{eqnarray}
where at the second inequality we have employed \Holder's inequality
and at the third we
have invoked the optimality of the coupling.  Another application of \Holder's
inequality, this time with conjugate exponents $\sfrac{p - 1}{p - 1 -
j}$ and $\sfrac{p - 1}{j}$, yields
\begin{equation}
\label{Hineq}
\| X^{p - 1 - j} Y^j \|_{p / (p - 1)} \leq \| X \|^{p - 1 - j}_p\,\| Y \|^j_p
\end{equation}
for $1 \leq j \leq p - 2$, and~\eqref{Hineq} is trivially an equality
when $j = 0$ or
$j = p - 1$.  Combining~\eqref{intermediate} and~\eqref{Hineq} and
rearranging, we
obtain the desired result.
\end{proof}

\begin{proof}[Proof of \refT{T:LBp}]
By \refL{L:decrease}, we have the bound
$$
\| Z_n\|_p \le d_p(Z_n,Y)+ \|Y\|_p \le d_p(Z_0,Y)+ \|Y\|_p
\le \|Z_0\|_p + 2\|Y\|_p.
$$
Thus Lemmas \refand{dpmoments}{L:persist} yield the
explicit bound
$$
d_p(Z_n,Y)\ge
\frac{| \EE\,Z^p_0 - \EE\,Y^p |}{\sum_{j = 0}^{p - 1} \left( \| Z_0
\|_p + 2 \| Y \|_p
\right)^{p - 1 - j}\,\| Y \|^j_p} \biggpar{\frac2{p+1}}^n .
\eqno\qedsymbol
$$
\noqed
\end{proof}

\begin{proof}[Proof of \refT{T:sharper}]
The 
cases $p = 1$ and $p = 2$ follow immediately from
\refT{T:LBp}; see also \refR{R:d12}.

When 
$p \geq 3$, fix $q>p$. By
Lemmas \refand{L:inter}{L:decrease} (the latter applied to~$d_q$),
for some $C_q$ we have
$$
d_p(F_n, F) \leq
C_q\, d_2^{\frac{2 (q - p)}{p (q - 2)}}(F_n, F)
$$
and thus \refT{T:LBp} implies \refT{T:sharper} with
$r=\bigpar{\frac{2}{p + 1}}^{\frac{p(q-2)}{2(q-p)}}$.
By taking $q$ sufficiently large, we obtain the result for any
$r<\bigpar{\frac{2}{p + 1}}^{p / 2}$.
\end{proof}

We have assumed in Theorems \refand{T:LB}{T:sharper} that
$F_0$ has finite moments of all orders. What happens if this fails?
If $\EE\,|Z_0|^p = \infty$ for some $p>2$,
then $d_p(Z_n,Y)=\infty$ for all $n$, but what can be said about
$d_2(Z_n,Y)$? It seems reasonable to conjecture that we have at least
as large $d_2$-distance in this case as in the nicer case with all
moments finite, and that Theorems \refand{T:LB}{T:sharper} hold
for all $F_0\neq F$.
 Unfortunately, we have not been able to prove this, but we offer 
the following partial result.

\begin{theorem}\label{T:anyr}
If $d_2(F_n,F)=O(r^n)$ for every $r>0$, then $F_0=F$.
Consequently, if $F_0\neq F$, there exists $r>0$ such that $d_2(F_n,F)
>r^n$ for infinitely many  values of $n$.
\end{theorem}
\begin{proof}
We may assume that
$\EE Z_0=0$ and $\EE Z_0^2=\EE Y^2$ 
(in particular, $\EE Z_0^2<\infty$),
because otherwise 
$d_2(F_n,F)=\Omega(r^n)$
with $r=2/3$:\ see \refR{R:d12}.  
By induction then 
$\EE Z_n=0$ and $\EE Z_n^2=\EE Y^2<1$ 
for every  $n$:\ see \refL{L:persist}.

As
usual, let $Z$, $Z^* \Leq Z$, and $U$ be independent.
If $|Z|\ge 2x$, $|Z^*|\le 2$, and $\frac23 \le U\le 1$, where $x\ge5$,
then
$$
|SZ|=|UZ+(1-U)Z^*+g(U)| \ge \tfrac23|Z|-\tfrac13|Z^*|-1
\ge \tfrac43 x-\tfrac53
\ge x.
$$
Thus,
$$
\PP(|SZ|\ge x) \ge \PP(|Z|\ge 2x)\cdot \PP(|Z|\le2)\cdot\tfrac 13,
\qquad x\ge5.
$$
If further $\EE Z^2\le1$, and thus by Chebyshev's inequality
$\PP(|Z|\le2)=1-\PP(|Z |>2)\ge 1 - \frac14 = \frac34$, this yields
$$
\PP(|SZ|\ge x) \ge \tfrac14\PP(|Z|\ge2x),\qquad x\ge5.
$$

Hence,
by induction on $n$ and our assumption on the first two moments
of $Z_0$, for any $x\ge5$,
$$
\PP(|Z_n|\ge x) \ge 4^{-n} \PP(|Z_0|\ge 2^nx), \qquad n\ge0,
$$
and in particular
\begin{equation}
\label{emma}
\PP(|Z_n|\ge 2^n) \ge 4^{-n} \PP(|Z_0|\ge 4^n), \qquad n\ge3.
\end{equation}

Now suppose that $d_2(Z_n,Y)=O(r^n)$. Using an optimal coupling
between $Z_n$ and $Y$, and the fact that $Y$ has moments of all
orders, we find
\begin{equation*}
\begin{split}
\PP(|Z_n|\ge 2^n) &\le
\PP(|Z_n-Y|\ge 2^{n-1})+\PP(|Y|\ge 2^{n-1}) \\
&\le 2^{2-2n} d_2^2(Z_n,Y)+\PP(|Y|\ge 2^{n-1})
=O(2^{-2n}r^{2n}).
\end{split}
\end{equation*}
Combining this with \eqref{emma}, we obtain, for $n\ge3$,
$$
\PP(|Z_0|\ge 4^n)
\le 4^n \PP(|Z_n|\ge 2^n)
=O(r^{2n}),
$$
which implies that $\EE|Z_0|^p<\infty$ for every $p>0$ such that $4^pr^2<1$.

Consequently,
if $d_2(Z_n,Y)=O(r^n)$
for every $r>0$, then $\EE|Z_0|^p<\infty$ for every $p>0$,
and \refT{T:sharper} applies to yield $F_0=F$.
\end{proof}

\begin{remark}
Our proof of \refT{T:anyr}, combined with the proof of \refT{T:sharper},
shows that if $F_0 \neq F$ and~$p$ (assumed $\geq 3$ here) is the smallest
positive integer such that either $\EE\,Z_0^p$ does not exist or
$\EE\,Z_0^p \neq \EE\,Y^p$, then
$d_2(Z_n,Y) > r^n$ for infinitely many values of $n$
for any $0 < r < r_p$, with
$r_p := 2^{-q}$, where $q$ is the unique solution in $(p, \infty)$ to
$2^\frac{2q(q-p)}{p(q-2)}=\frac{p+1}2$.
\end{remark}

\section{Other lower bounds}\label{S:otherlower}

In \refS{S:lower} we showed that convergence of the iterates~$F_n$ to~$F$
in the $d_2$-metric is not faster than geometric.
In this final section we show likewise that
the convergence is not faster than geometric
in the other metrics we have considered in this paper.
We again assume that $F_0 \neq F$ has finite moments of all orders.
(Without this
hypothesis, we can prove partial results by the method used in the
proof of \refT{T:anyr}, but we do not know whether the full results hold.)

\subsection{\KSm{} and total variation distances}

We begin with a simple lemma.

\begl\label{L:ksmoments}
Let $p>0$.
For any $X \sim F$ and $Y \sim G$ each with finite {\rm $p$th}
absolute
moment,
if $K :=\dks(F, G)$, then, for any $0 \leq M < \infty$,
$$
\left| \EE \left( X^p; X > 0 \right) - \EE \left( Y^p; Y > 0 \right)
\right| \leq  K M^p
+  \EE \left( X^p; X > M \right)+ \EE \left( Y^p; Y > M\right)
$$
and, if $p$ is an integer,
$$
\left| \EE\,X^p - \EE\,Y^p \right| \leq 2 K M^p
+ \EE \left( |X|^p; |X| > M \right) + \EE \left( |Y|^p; |Y| > M \right).
$$
\enl

\begin{proof}
Define $X_M:=\min(X^+,M)$, where $X^+=\max(X,0)$, and similarly $Y_M$.
Then
$$
0 \le
\EE(X^p; X>0) - \EE X_M^p
\le \EE(X^p; X>M)
$$
and similarly for $Y$, while
\begin{equation*}
\begin{split}
\left| \EE\,X_M^p - \EE\,Y_M^p \right|
&= \left|\int_0^M p x^{p-1} \PP(X>x)\, dx
  -\int_0^M p x^{p-1} \PP(Y>x)\, dx \right|\\
&\le\int_0^M p x^{p-1} \bigl|\PP(X>x)-\PP(Y>x)\bigr|\, dx
\le KM^p.
\end{split}
\end{equation*}

Together, these yield the first inequality.

The second follows by applying the first to $(X,Y)$ and to $(-X,-Y)$ and
summing or subtracting, depending on the parity of $p$.
\end{proof}

\begt\label{T:kslower}
Suppose $F_0 \neq F$ has finite moments of all orders,
and let~$p$ be
defined as in \refT{T:sharper}.
Then, for any $0< r < 2 / (p + 1)$,
$$
\dtv(F_n,F)\ge
\dks(F_n, F) = \Omega(r^n).
$$
(The implicit multiplicative constant depends on
both~$F_0$ and the choice of~$r$.)
\ent

\begin{proof}
Let $K_n := \dks(F_n, F) > 0$.
If we apply
\refL{L:ksmoments}
and then use \refL{L:persist}, we find, for any $q\ge p$,
\begin{equation}
\label{l1}
\begin{split}
\left| \EE\,Z^p_0 - \EE\,Y^p \right| \bigpar{\sfrac{2}{p + 1}}^n
&\leq 2 K_n M^p
 + \EE \left( |Z_n|^p; |Z_n| > M \right) + \EE \left( |Y|^p; |Y| > M \right)\\
&\leq 2 K_n M^p
 + M^{-(q-p)}\EE |Z_n|^q  + M^{-(q-p)}\EE |Y|^q,
\end{split}
\end{equation}
for any $0 \leq M < \infty$.
It follows from \refL{L:decrease} that $\EE |Z_n|^q\le C_q$, for some
$C_q$ not depending on $n$. Choosing $M=K_n^{-1/q}$ thus gives, with
$c := \left| \EE\,Z^p_0 - \EE\,Y^p \right|>0$,
\begin{equation*}
c\bigpar{\sfrac{2}{p + 1}}^n
\leq 2 K_n^{1 - (p/q)} + 2C_q K_n^{1 - (p/q)},
\end{equation*}
and
thus $K_n=\Omega(r^n)$ with
$r=\bigpar{\frac{2}{p + 1}}^{q/(q-p)}$.
The result follows, since $r\to 2/(p+1)$ as $q\to\infty$.
\end{proof}

\subsection{Density functions and characteristic functions}
We immediately obtain results for the density functions $f_n$, which
by \refT{T:density} exist at least for $n\ge3$ and by
\refT{T:densconv} converge uniformly, at a geometric rate,
to the density $f$ of $Y$.

\begin{corollary}\label{C:denslower}
Suppose $F_0 \neq F$ has finite moments of all orders,
and let~$p$ be
defined as in \refT{T:sharper}.
Then, for any $0< r < 2 / (p + 1)$,
\begin{equation}\label{dl1}
\int_{-\infty}^\infty|f_n(x)-f(x)|\,dx =\Omega(r^n)
\end{equation}
and
\begin{equation}\label{dl2}
\supfnf = \Omega(r^n).
\end{equation}
\end{corollary}

\begin{proof}
The estimate \eqref{dl1} follows from \refT{T:kslower}, because
(whenever $f_n$ exists)
$\int_{-\infty}^\infty|f_n(x)-f(x)|\,dx =2\dtv(F_n,F)$.

The estimate \eqref{dl2} follows from \refT{T:kslower} using
inequality~\eqref{g1} and the discussion following it.
\end{proof}

Similarly, we have a geometric lower bound for the  $L^1(\RR)$ and
$L^\infty(\RR)$
distances of the characteristic functions.

\begin{corollary}\label{C:charlower}
Suppose $F_0 \neq F$ has finite moments of all orders,
and let~$p$ be
defined as in \refT{T:sharper}.
Then, for any $0< r < 2 / (p + 1)$,
\begin{equation}\label{cl1}
\int_{-\infty}^\infty|\phi_{Z_n}(t)-\phi_Y(t)|\,dt =\Omega(r^n).
\end{equation}
and
\begin{equation}\label{cl2}
\sup_{t}|\phi_{Z_n}(t)-\phi_Y(t)| =\Omega(r^n).
\end{equation}
\end{corollary}
\begin{proof}
The estimate \eqref{cl1} is immediate from \refC{C:denslower} and
inequality~\eqref{diff1}.
Next, \eqref{cl1} for some $r=r_0$ and \refT{T:phi} imply
\eqref{cl2} for any $r<r_0$ by the argument used to show
\eqref{axel} in the proof of \refT{T:densconv}.
\end{proof}

It is not too hard to extend Corollaries \refand{C:denslower}{C:charlower}
to any $L^q(\RR)$ distance, $1 \le q \le \infty$.

\subsection{Moment generating functions}
Finally, we consider 
lower bounds for the convergence of moment
generating functions.
We assume for simplicity that $Z_0$ has an everywhere finite \mgf, and
know by \refT{T:mgf}
that then $\psi_{Z_n}$ converges to $\psi_Y$
pointwise, and uniformly on compact sets, with geometric rate.
For lower bounds we first note that
if $F_0\neq F$ and $p$ is as in \refT{T:sharper}, then
the derivatives of $\psi_{Z_n}$ and $\psi_Y$ at the origin,
which equal the corresponding moments of $Z_n$ and $Y$,
by \refL{L:persist} coincide
up to order $p-1$, while the $p$th derivatives differ by
$c\bigpar{\frac2{p+1}}^n$ with $c=\EE Z_0^p-\EE Y^p\neq0$.
For $\gl$ close to the origin, this and a Taylor expansion
shows that
$|\psi_{Z_n}(\gl)-\psi_Y(\gl)|=\Omega\bigpar{|\gl|^p\bigpar{\frac2{p+1}}^n}$,
but the range of $\gl$ where we can prove that this is valid
depends on $n$.

Indeed, there is no general lower bound for
$|\psi_{Z_n}(\gl)-\psi_Y(\gl)|$
for a fixed $\gl$, since
there may be points~$\gl$ where $\psi_{Z_n}(\gl)$ and $\psi_Y(\gl)$
coincide ``accidentally''.
For example, suppose that $Z_0$ is bounded with $\EE\,Z_0 = 0$ and
$\Var\,Z_0 > \Var\,Y$.
By induction, the same holds for
each~$Z_n$: see \refL{L:persist}
and note that $g(U)$ is bounded.
Consequently, for each $n$, Taylor's formula shows that
$\psi_{Z_n}(\gl)>\psi_Y(\gl)$ for small positive~$\gl$,
while $\psi_{Z_n}(\gl)=\exp\bigpar{O(\gl)}$ and thus \refT{T:mgf_lower} shows
that $\psi_{Z_n}(\gl) < \psi_Y(\gl)$ for large $\gl$. Hence there exists for
every $n$ at least one positive $\gl = \gl_n$ such that
$\psi_{Z_n}(\gl) = \psi_Y(\gl)$.
Nevertheless, such points have to be isolated, and if we consider the
maximum deviation over an interval, we
have a geometric lower bound.

\begt\label{mgflower}
Suppose $F_0 \neq F$ has everywhere finite moment generating function, and let
$(a, b)$ be a nonempty
interval. Then there exists $r>0$ such that
$\sup_{a \le \gl \le b}|\psi_{Z_n}(\gl)-\psi_Y(\gl)| = \Omega(r^n)$.
\ent

\begin{proof}
We
use the fact that the \mgf{s} $\psi_{Z_n}$ and $\psi_Y$ are entire
analytic functions in the complex plane $\CC$.

Let $R := |a| + |b| + 1$. There exists a (unique) function $\omega$ which is
continuous on $ D_R:=\set{z\in\CC:|z|\le R}$ and analytic in
$\Omega := \set{z:|z|<R}\setminus[a,b]$ such that $\omega(z)=0$ for
$|z|=R$ and $\omega(z)=1$ for $z\in[a,b]$; this function is called
harmonic measure and is probabilistically given by the probability
that a Brownian motion starting at $z$ hits $[a,b]$ before it hits
$\set{z:|z|=R}$.

Let $f_n(z) := \psi_{Z_n}(z) - \psi_Y(z)$
and $u_n(z) := \ln |f_n(z)| \ge -\infty$.
For $z\in  D_R$,
$$
|f_n(z)|\le |\psi_{Z_n}(z)|+|\psi_Y(z)|
\le \psi_{Z_n}(R) + \psi_Y(R) + \psi_{Z_n}(-R) + \psi_Y(-R),
$$
which by \refT{T:mgf} is bounded by some constant $A<\infty$
(depending on $Z_0$ but not on $n$).
Let further $\gd_n := \max_{a \le \gl \le b} |f_n(\gl)|$;
we may of course restrict attention to those values of~$n$
satisfying $\gd_n < 1$.
Now $u_n(z)\le \ln A$ for $|z|=R$ and $u_n(z) \le \ln \gd_n$ for
$z\in [a,b]$; thus (since $A\ge1$)
\begin{equation}
\label{major}
u_n(z) \le \ln A + (\ln \gd_n) \omega(z)
\end{equation}
for every $z\in \partial \Omega$.  Since $u_n$ is subharmonic and
the right hand side is harmonic in $\Omega$ and continuous on its
closure, \eqref{major} holds for every
$z\in \overline\Omega=D_R$,
\cf{} \cite[Theorems 17.3 and 17.4]{RudinRealComplex}.
In particular, setting $\eps:=\inf_{|z|\le1} \omega(z)>0$, we have
$$
u_n(z) \le \ln A + \eps \ln \gd_n , \qquad |z| \le 1,
$$
or
\begin{equation}
\label{fmax}
|f_n(z)| \le  A \gd_n^\eps , \qquad |z| \le 1.
\end{equation}
Let $p$ be as in \refT{T:sharper}. By \eqref{fmax} and Cauchy's
estimates
\cite[Theorem 10.26]{RudinRealComplex},
$$
|f_n^{(p)}(0)| \le p!\,A \gd_n^\eps.
$$
Since by \refL{L:persist}
$$
|f_n^{(p)}(0)|=|\EE Z_n^p-\EE Y^p|=\Omega\bigpar{(\tfrac{2}{p + 1})^n},
$$
it follows that $\gd_n = \Omega(r^n)$ with $r=\bigpar{\frac2{p+1}}^{1/\eps}$.
\end{proof}

\end{document}